\documentclass[reqno,centertags, 12pt]{amsart}
\usepackage{amsmath,amsthm,amscd,amssymb}
\usepackage{latexsym}
\usepackage{graphicx}

\newcommand{\bbC}{{\mathbb{C}}}
\newcommand{\bbD}{{\mathbb{D}}}
\newcommand{\bbE}{{\mathbb{E}}}
\newcommand{\bbL}{{\mathbb{L}}}

\newcommand{\bbR}{{\mathbb{R}}}
\newcommand{\bbS}{{\mathbb{S}}}
\newcommand{\bbU}{{\mathbb{U}}}
\newcommand{\bbZ}{{\mathbb{Z}}}
\newcommand{\bdone}{{\boldsymbol{1}}}
\newcommand{\bddot}{{\boldsymbol{\cdot}}}

\newcommand{\dott}{\,\cdot\,}

\newcommand{\lb}{\label}
\newcommand{\f}{\frac}

\newcommand{\ol}{\overline}
\newcommand{\ti}{\tilde  }

\newcommand{\tr}{\text{\rm{Tr}}}
\newcommand{\dist}{\text{\rm{dist}}}
\newcommand{\intt}{\text{\rm{int}}}

\newcommand{\spec}{\text{\rm{spec}}}

\newcommand{\ess}{\text{\rm{ess}}}

\newcommand{\s}{\text{\rm{s}}}

\newcommand{\supp}{\text{\rm{supp}}}

\newcommand{\HS}{\text{\rm{HS}}}
\newcommand{\bi}{\bibitem}

\newcommand{\beq}{\begin{equation}}
\newcommand{\eeq}{\end{equation}}
\newcommand{\ba}{\begin{align}}
\newcommand{\ea}{\end{align}}
\newcommand{\veps}{\varepsilon}

 \DeclareMathOperator{\Ima}{Im}

 \allowdisplaybreaks
\numberwithin{equation}{section} \DeclareMathOperator*{\Sz}{Sz}

\newtheorem{theorem}{Theorem}[section]

\newtheorem{proposition}[theorem]{Proposition}
\newtheorem{lemma}[theorem]{Lemma}
\newtheorem{corollary}[theorem]{Corollary}
\theoremstyle{definition}
\newtheorem{example}[theorem]{Example}

\newtheorem{oq}[theorem]{Open Question}

\theoremstyle{remark}
\newtheorem*{remark}{Remark}
\newtheorem*{remarks}{Remarks}
\theoremstyle{definition}
\newtheorem*{definition}{Definition}
\newcommand{\abs}[1]{\lvert#1\rvert}

\newcounter{smalllist}
\newenvironment{SL}{\begin{list}{{\rm\roman{smalllist})}}{%
\setlength{\topsep}{0mm}\setlength{\parsep}{0mm}\setlength{\itemsep}{0mm}%
\setlength{\labelwidth}{2em}\setlength{\leftmargin}{2em}\usecounter{smalllist}%
}}{\end{list}}

\begin{document}

\title[Fine Structure of the Zeros of OP, IV]
{Fine Structure of the Zeros of Orthogonal Polynomials, \\
IV.\ A Priori Bounds and Clock Behavior}
\author[Y.~Last and B.~Simon]{Yoram Last$^{1,3}$ and Barry Simon$^{2,3}$}

\thanks{$^1$ Institute of Mathematics, The Hebrew University,
91904 Jerusalem, Israel. E-mail: ylast@math.huji.ac.il. Supported in part
by The Israel Science Foundation (Grant No.\ 188/02)}
\thanks{$^2$ Mathematics 253-37, California Institute of Technology, Pasadena, CA 91125.
E-mail: bsimon@caltech.edu. Supported in part by NSF grant DMS-0140592}
\thanks{$^3$ Research supported in part
by Grant No.\ 2002068 from the United States-Israel Binational Science Foundation
(BSF), Jerusalem, Israel}

\dedicatory{Dedicated to Percy Deift on the occasion of his sixtieth birthday}

\date{May 16, 2006}

\begin{abstract} We prove locally uniform spacing for the zeros of orthogonal
polynomials on the real line under weak conditions (Jacobi parameters approach
the free ones and are of bounded variation). We prove that for ergodic discrete
Schr\"odinger operators, Poisson behavior implies positive Lyapunov exponent. Both
results depend on a priori bounds on eigenvalue spacings for which we provide
several proofs.
\end{abstract}

\maketitle

\section{Introduction} \lb{s1}

Our primary goal in this paper concerns the fine structure of the zeros of
orthogonal polynomials on the real line (OPRL), although we will say something
about zeros of paraorthogonal polynomials on the unit circle (POPUC) (see
Section~\ref{s10}). Specifically, $d\mu$ will be a measure of compact support which
is nontrivial (i.e., not supported on a finite set), usually a probability measure.
$P_n(x)$ or $P_n(x;d\mu)$ will be the monic orthogonal polynomials and $p_n(x)=
P_n/\|P_n\|$ the orthonormal polynomials. The Jacobi parameters, $\{a_n,b_n\}_{n=1}^\infty$,
are defined by
\begin{equation} \lb{1.1}
xP_n(x)=P_{n+1}(x) + b_{n+1} P_n(x) + a_n^2 P_{n-1}(x)
\end{equation}
$n=0,1,2,\dots$, where $P_{-1}(x)\equiv 0$. It follows that (when $\mu(\bbR)=1$)
\begin{equation} \lb{1.2}
\|P_n\|=a_1\dots a_n
\end{equation}
so
\begin{equation} \lb{1.3}
xp_n(x)=a_{n+1} p_{n+1}(x) + b_{n+1} p_n(x) + a_n p_{n-1} (x)
\end{equation}
and $x$ has a matrix representation in the orthonormal basis $\{p_n\}_{n=0}^\infty$,
\begin{equation} \lb{1.4}
J=\begin{pmatrix} b_1 & a_1 & 0 & \cdots \\
a_1 & b_2 & a_2 & \cdots \\
\vdots & \vdots & \ddots & \ddots
\end{pmatrix}
\end{equation}
called the Jacobi matrix. The finite Jacobi matrix, $J_{n;F}$, is the $n\times n$
submatrix of $J$ in the upper left corner. It is easy to see that (see (b) of the
Appendix)
\begin{equation} \lb{1.5}
P_n(x) =\det (x\bdone -J_{n;F})
\end{equation}

Let $\{x_j^{(n)}\}_{j=1}^n$ be the zeros of $P_n(x)$, so \eqref{1.5} says that the
$x_j^{(n)}$ are eigenvalues of $J_{n;F}$. Let $d\nu_n$ be the pure point probability
measure that gives weight $1/n$ to each $x_j^{(n)}$. We say the {\it density of states\/}
exists if $d\nu_n$ has a weak limit $d\nu_\infty$. By \eqref{1.5}, one sees that
\begin{equation} \lb{1.5a}
\int x^k\, d\nu_n(x) = \f{1}{n}\, \tr (J_{n;F}^k)
\end{equation}
which is often useful.

The existence of the limit for a large class of regular measures on $[-2,2]$ goes back
to Erd\"os-Tur\'an \cite{ET40}. Nevai \cite{Nev79} realized all that was used is
\begin{equation} \lb{1.6}
\lim_{n\to\infty}\, \abs{b_n} + \abs{a_n-1} =0
\end{equation}
(following the convention in the Schr\"odinger operator community, we use $a_n\equiv 1$
as a ``free" case, while the OP community uses $a_n\equiv \f12$). Indeed,

\begin{theorem}[known]\lb{T1.1} If $d\mu$ is a measure where \eqref{1.6} holds,
the density of states exists and is given by
\begin{equation} \lb{1.7}
d\nu = \pi^{-1} \bigl(\sqrt{4-x^2}\,\bigr)^{-1} \chi_{[-2,2]}\, dx
\end{equation}
\end{theorem}

The modern proof is not hard. For $a_n\equiv 1$, $b_n\equiv 0$, the OPRL are explicitly
given by
\[
P_n (2\cos\theta) = \f{\sin ((n+1)\theta)}{\sin\theta}
\]
from which one computes $d\nu$ exactly for this case. If $J_{0,n;F}$ is the corresponding
cutoff $J_0$, then \eqref{1.6} implies
\begin{equation} \lb{1.8}
\f{1}{n}\, \tr (J_{n;F}^k - J_{0,n;F}^k)\to 0
\end{equation}
for each $k=0,1,2,\dots$ which, by \eqref{1.5a}, implies $\int x^k\, d\nu_n$ has the same
limits as for the free case.

Another case where it is known that $d\nu$ exists is when $a_n,b_n$ are samples of an
ergodic family, that is, $a_n^{(\omega)}=f(T^n \omega)$, $b_n^{(\omega)}=g(T^n \omega)$
with $T\colon \Omega\to\Omega$ an ergodic transformation on $(\Omega,d\rho)$, a probability
measure space. In that case, it is known (going back to the physics literature and proven
rigorously by Pastur \cite{P73}, Avron-Simon \cite{S149}, and Kirsch-Martinelli
\cite{KM}):

\begin{theorem}[known]\lb{T1.2} For ergodic Jacobi matrices, $d\nu_{n,\omega}$ has a limit
$d\nu$ for a.e.\ $\omega$ and $d\nu$ is a.e.\ $\omega$-independent.
\end{theorem}

Again, the proof uses \eqref{1.5a} plus in this case that, by ergodicity, $\f{1}{n}
\tr (J_{n;F}^k)$ has a limit a.e.\ by the Birkhoff ergodic theorem.

The most important examples of the ergodic case are periodic, almost periodic,
and random.

One easily combines the two ideas to see that $d\nu$ exists (and does not depend on
$\delta a_n, \delta b_n$) if
\[
a_n = a_n^{(\omega)} + \delta a_n \qquad
b_n = b_n^{(\omega)} + \delta b_n
\]
with $a_n^{(\omega)}, b_n^{(\omega)}$ ergodic  and $\abs{\delta a_n} + \abs{\delta b_n}
\to 0$.

These results describe the bulk features of the zeros. Here we are interested in the
fine structure,  on the level of individual eigenvalues; specifically, the focus in
\cite{Saff1,Saff2,Saff3} and a main focus in this paper is what we call clock
behavior, that the spacing locally is equal spacing. The term clock comes from the
case of orthogonal polynomials on the unit circle (OPUC) where $d\nu$ is typically
Lebesgue measure on a circle and the equal space means the zeros look like the
numbers on a clock.

In order for the density of zeros to be $d\nu$, the equal spacing must be $1/(d\nu/dE)$.
The symmetric derivative $d\nu/dE$ exists for a.e.\ $E$ and, of course, $(d\nu/dE)\,dE$
is the a.c.\ part of $d\nu$. It is known (see, e.g., Avron-Simon \cite{S149}) that $d\nu$
has no pure points and, in many cases, it is known that $d\nu/dE$ is a continuous function,
or even $C^\infty$. To be formal, we define first

\begin{definition} Let $E_0\in\supp(d\nu)$. We let $z_n^{(j)}(E_0)$ be the zeros nearest
$E_0$ so that
\[
z_n^{(-2)} (E_0) < z_n^{(-1)}(E_0) \leq E_0 < z_n^{(1)}(E_0) < z_n^{(2)}(E_0) < \cdots
\]
if such zeros exist. If $E_0\in[\supp(d\nu)]^{\intt}$, then $z_n^j (E_0)$ exists for
$j$ fixed and $n$ large.
\end{definition}

\begin{definition} Let $E_0\in \supp(d\nu)$. We say there is weak clock behavior at $E_0$
if $d\nu/dE_0$ exists and
\begin{equation} \lb{1.9}
\lim_{n\to\infty}\, n [z_n^{(1)}(E_0) - z_n^{(-1)}(E_0)]\, \f{d\nu}{dE}\, (E_0) =1
\end{equation}
We say there is strong clock behavior at $E_0$ if $d\nu/dE_0$ exists, \eqref{1.9} holds,
and for $j=1,\pm 2, \pm 3, \dots$ fixed,
\begin{equation} \lb{1.10}
\lim_{n\to\infty}\, n [z_n^{(j)}(E_0) - z_n^{(j+1)} (E_0)] \, \f{d\nu}{dE}\, (E_0) =1
\end{equation}
\end{definition}

\begin{definition} We say there is uniform clock behavior on $[\alpha,\beta]$ if $d\nu/dE$
is continuous and nonvanishing on $[\alpha,\beta]$ and
\begin{equation} \lb{1.11}
\lim_n \biggl[\sup \biggl\{ \biggl| n[E-E']-\biggl(\f{d\nu}{dE}\biggr)^{-1}\biggr|
\biggm| E,E' \text{ are successive zeros of $p_n$ in } [\alpha,\beta]\biggr\}\biggr] = 0
\end{equation}

It is obvious that uniform clock behavior implies strong clock behavior at each interior
point.
\end{definition}

In the earlier papers in this series that discussed clock behavior for OPRL
\cite{Saff1,Saff3}, there was a technical issue that severely limited the results
in general situations. In all cases, a Jost function-type analysis was used to show
that a suitably rescaled $p_n$ converged, that is,
\begin{equation} \lb{1.13a}
c_n p_n \biggl( E_0 + n(x-E_0) \, \f{d\nu}{dE}\, (E_0) + \xi_n\biggr) \to f_\infty (x)
\end{equation}
where $\abs{\xi_n}\leq c/n$ and $f_\infty$ has zeros at $0,\pm 1, \pm 2, \dots$. This would
naively seem to show that $p_n$ has clock-spaced zeros and, indeed, it does imply at least
one zero near $E_0 + \xi_n + \f{j}{n} (\f{d\nu}{dE_0})^{-1}$ consistent with clock spacing.

The snag involves uniqueness, for the function
\begin{equation} \lb{1.12}
\sin (\pi x) - \f{1}{n} \sin (n^2 \pi x)
\end{equation}
has a limit like $f_\infty$ but has more and more zeros near $x=0$. That is, one needs to
prove uniqueness of the zeros near $E_0 + \xi_n + \f{j}{n} (\f{d\nu}{dE_0})^{-1}$.

In previous papers in this series, two methods were used to solve this uniqueness problem.
One relied on some version of the argument principle, essentially Rouch\'e's theorem. This
requires analyticity which, typically, severely restricts what recursion coefficients
are allowed. In the case of OPUC where one needs to control zeros in the complex plane,
some kind of analyticity argument seems to be necessary. The second method relies on the
fact that if \eqref{1.13a} also holds for derivatives and $f'_\infty (j)\neq 0$, then
there is a unique zero. This argument also requires extra restrictions on the recursion
coefficients, albeit not as severe as the analyticity requirement. For example, in
\cite{Saff1}, one only needed $\sum_{n=1}^\infty \abs{a_n-1} + \abs{b_n} <\infty$
to get \eqref{1.13a} for $E\in [-2+\veps, 2-\veps]$. But to control derivatives, we
needed $\sum_{n=1}^\infty n(\abs{a_n-1} + \abs{b_n})<\infty$. In fact, the general
argument failed to capture Jacobi polynomials (whose clock estimates were earlier
obtained by V\'ertesi and Szabados \cite{SV,Ver90,Ver97,Ver99}) where separate
arguments are needed. In \cite{Saff3}, it was decided not to deal with asymptotically
periodic OPRL since the derivative arguments looked to be tedious.

\smallskip
{\it Note.} The papers \cite{Saff1,Saff2,S303} discussed asymptotics of OPUC when
the Verblunsky coefficients decay exponentially. We have discovered a paper of
Pan \cite{Pan96} whose results overlap those in \cite{Saff1,Saff2,S303}.

\smallskip
The key realization of this paper is that there is a more efficient way to eliminate
pathologies like those in \eqref{1.12}. Namely, we will seek a priori lower bounds on
eigenvalue spacings. If we find any $O(1/n)$ lower bound, that implies the rescaled $p_n$
of \eqref{1.13a} has at most one zero near any zero of $f_\infty$. Such lower bounds are
not new in suitable situations. Erd\"os-Tur\'an \cite{ET40} already have such bounds if
the measure is purely absolutely continuous in an interval with a.c.\ weights bounded
away from zero and infinity. These ideas were developed by Nevai \cite{Nev79} and Golinskii
\cite{Gol02}. Under suitable hypotheses on the transfer matrix, lower bounds are known in
the Schr\"odinger operator community; see, for example, Jitomirskaya et al.\ \cite{JSS03}.
({\it Note}: They use ``equal spacing" for $O(1/n)$ lower bounds and do not mean clock
behavior by this term.)

While we could have used these existing bounds in the proofs of Theorems~\ref{T1.3}
and \ref{T1.4} below, we have found a new approach which allows us to also prove
Theorem~\ref{T1.5} below, and this approach is discussed in Section~\ref{s2}.
With these lower bound ideas, we can prove:

\begin{theorem}[$\equiv$ Theorem~\ref{T3.2}] \lb{T1.3} Suppose that
\begin{equation} \lb{1.14}
\sum_{n=1}^\infty (\abs{a_n-1} + \abs{b_n}) <\infty
\end{equation}
Then there is uniform clock behavior on each interval $[-2+\veps, 2-\veps]$ for
any $\veps >0$.
\end{theorem}

\begin{theorem}[$\equiv$ Theorem~\ref{T3.4}] \lb{T1.4} Suppose that
\begin{gather}
\lim_{n\to\infty}\, a_n =1 \qquad \lim_{n\to\infty}\, b_n =0 \lb{1.15} \\
\sum_{n=1}^\infty (\abs{a_{n+1} -a_n} + \abs{b_{n+1} -b_n}) <\infty \lb{1.16}
\end{gather}
Then there is uniform clock behavior on each interval $[-2+\veps, 2-\veps]$ for
any $\veps >0$.
\end{theorem}

\begin{remarks} 1. Theorem~\ref{T1.4} implies Theorem~\ref{T1.3}, but we state them
as separate theorems since the proof of the first is easier.

2. We will also prove results of this genre for perturbations of periodic recursion
coefficients and for where $n+1$ in \eqref{1.16} is replaced by $n+p$ for some $p$.

3. We also obtain results (see Theorem~\ref{T3.8}) near $\pm 2$ if
\begin{equation} \lb{1.17}
\sum_{n=1}^\infty n (\abs{a_n-1} + \abs{b_n}) <\infty
\end{equation}
\end{remarks}

Using our strong lower bound, we will also prove the following:

\begin{theorem} \lb{T1.5} Let $a_n^{(\omega)},b_n^{(\omega)}$ be ergodic Jacobi
parameters. Let $E_0$ be such that
\begin{SL}
\item[{\rm{(i)}}] The Lyapunov exponent $\gamma (E_0)=0$.
\item[{\rm{(ii)}}] The symmetric derivative of $\nu$ exists at $E_0$ and is finite
and nonzero.
\end{SL}
Then there exists $C>0$ so that with probability $1$,
\begin{equation} \lb{1.18}
\liminf_{n\to\infty}\, n [z_n^{(1)}(E_0) -z_n^{(-1)}(E_0)] >C
\end{equation}
\end{theorem}

This result is especially interesting because it implies that the zeros cannot have
Poisson behavior. It implies that Poisson behavior {\it and ergodicity\/} require
$\gamma (E_0) >0$. We will say more about these issues in Sections~\ref{s4} and \ref{s11}.

While our initial interest in a priori lower bounds came from clock theorems
and we could have finished the paper with Section~\ref{s4} if our sole purpose
was to prove Theorems~\ref{T1.3}--\ref{T1.5}, it seemed natural to also consider
upper bounds. Moreover, in looking over the upper/lower bound results in the OP
literature, we realized one could get more from these methods, so we discuss
that also.

Broadly speaking, we have two sets of results and methods. The methods rely on
either transfer matrices with hypotheses on recursion coefficients or on OP methods
with hypotheses on the measure. We believe that the results are of interest to
both the Schr\"odinger operator and OP communities. Because the OP methods
are unfamiliar to many Schr\"odinger operator theorists and are easy to prove
(albeit very powerful), we have included an appendix with some major OP methods.

The detailed plan of the paper is as follows: In Section~\ref{s2}, we prove our
a priori lower bounds involving the transfer matrix (or, more precisely, the
growth of subordinate and nonsubordinate solutions). In Section~\ref{s3}, we
prove a variety of clock theorems, including Theorems~\ref{T1.3} and \ref{T1.4}.
In Section~\ref{s4}, using ideas of Deift-Simon \cite{S169}, we prove
Theorem~\ref{T1.5}. In Section~\ref{s5}, we obtain upper bounds on eigenvalue
spacing using the transfer matrix. Section~\ref{s6} discusses using suitable Pr\"ufer
angles to control spacing of zeros.

In Section~\ref{s7}, we begin our discussion of OP methods with a technical
result on $L^p$ bounds on $w^{-1}$ and bounds of the Christoffel function.
These bounds, which we will need for examples later are local versions of
some bounds of Geronimus \cite{GBk} with a rather different method of proof.
In Section~\ref{s8}, we discuss upper bounds on eigenvalue spacing using OP
methods and, in particular, find a remarkable lower bound on the density of
states that is a kind of microlocal version of some bounds of Deift-Simon
\cite{S169}. In Section~\ref{s9}, we discuss lower bounds on eigenvalue spacing.
The methods in Sections~\ref{s8} and \ref{s9} are borrowed from Erd\"os-Tur\'an
\cite{ET40,Szb}, Nevai \cite{Nev79}, and Golinskii \cite{Gol02}, but we show how
to localize them and squeeze out stronger results. In Section~\ref{s10}, we briefly
discuss the analogs of our results for zeros of POPUC, and in Section~\ref{s11}
discuss a number of examples, counterexamples, conjectures, and questions.

\smallskip
It is a pleasure to thank L.~Golinskii, S.~Jitomirskaya, R.~Killip, P.~Nevai,
and M.~Stoiciu for useful discussions. This research was begun during B.~S.'s
stay as a Lady Davis Visiting Professor at The Hebrew University of Jerusalem.
He would like to thank H.~Farkas for the hospitality of the Einstein Institute
of Mathematics at The Hebrew University, and E.~B.~Davies and A.~N.~Pressley
for the hospitality of King's College, London where this was work was continued.
Y.~L. would like to thank G.~A.~Lorden and T.~A.~Tombrello for the hospitality
of Caltech, where this work was completed.

Percy Deift has long been a player in spectral theory and more recently, a
champion for orthogonal polynomials. In particular, this paper exploits the
work of Deift-Simon \cite{S169}. It is a pleasure to dedicate this paper
to Percy.

\section{Variation of Parameters and Lower Bounds
via Transfer Matrices} \lb{s2}

Our goal here is to use variation of parameters to study eigenvalue spacing. Variation
of parameters has an ancient history going back to Lagrange \cite{HW} and it was
extensively used to study variation of solution with change in potential, for example,
to study asymptotics in tunnelling problems \cite{S105}. The usefulness of the method as
a tool in spectral theory
goes back at least to the work of Gilbert-Pearson \cite{GP}
(also see \cite{Gil,KP})
with significant later contributions by Jitomirskaya-Last \cite{JL0,JL1,JL2}
and Killip-Kiselev-Last \cite{KKL}. It is essentially their equation we will
use although, interestingly enough, the earlier applications are to $E,E'$
with $E\in\bbR$ and $E'=E+i\veps$, while our application is to $E,E'$ both
in $\bbR$.

Given $E\in\bbC$, we consider solutions of
\begin{equation} \lb{2.1}
a_n u_{n+1} + (b_n-E)u_n + a_{n-1} u_{n-1} =0
\end{equation}
for $n=1,2,\dots$. Here $\{a_n,b_n\}_{n=1}^\infty$ are the Jacobi parameters of the
measure we are considering and
\begin{equation} \lb{2.2}
a_0\equiv 1
\end{equation}

For $\theta\in [0,\pi)$, we denote by $u_n (E,\theta)$ the solution of \eqref{2.1} with
\begin{equation} \lb{2.3}
u_0 (E,\theta) =\sin(\theta) \qquad
u_1 (E,\theta) =\cos(\theta)
\end{equation}
In particular,
\begin{equation} \lb{2.4}
u_n (E,\theta=0)=p_{n-1}(E) \qquad n=1,2,\dots
\end{equation}
The transfer matrix is defined by
\begin{equation} \lb{2.5}
T(n,E) = \begin{pmatrix}
u_{n+1} (E,\theta=0) & u_{n+1} (E,\theta=\pi/2) \\
u_n (E,\theta=0) & u_n (E,\theta=\pi/2)
\end{pmatrix}
\end{equation}
so for any solution of \eqref{2.1},
\begin{equation} \lb{2.6}
\binom{u_{n+1}}{u_n} = T(n,E) \binom{u_1}{u_0}
\end{equation}

Let $K(n,m;E)$ be the kernel
\begin{equation} \lb{2.7}
K(n,m;E) = u_n (E,0) u_m (E, \pi/2) - u_n (E, \pi/2) u_m (E,0)
\end{equation}
Define the operator $A_L(E)$ on $\bbC^L=\{v_\bddot = \{v_k\}_{k=1}^L\mid v_k\in\bbC\}$ by
\begin{equation} \lb{2.8}
(A_L(E)v)_n =\sum_{m=1}^n K(n,m;E) v_m
\end{equation}
(note that $K(n,n)=0$, so the sum also goes to $n-1$). The following summarizes results
from \cite{JL0,JL1,JL2,KKL}:

\begin{theorem}\lb{T2.1} Let $w_n (E), w_n(E')$ solve \eqref{2.1} for $E,E'$, and
suppose
\begin{equation} \lb{2.9}
w_0 (E)=w_0 (E') \qquad w_1 (E)=w_1(E')
\end{equation}
Then
\begin{equation} \lb{2.10}
w_n (E') = w_n (E) + (E'-E) \sum_{m=1}^n K(n,m;E) w_m (E')
\end{equation}
that is,
\begin{equation} \lb{2.11}
w_\bddot (E') = w_\bddot (E) + (E'-E) (A_L(E) w(E'))_\bddot
\end{equation}
Moreover, with $\|\cdot\|_{\HS}$ the Hilbert-Schmidt norm on $\bbC^L$ and
$\|\cdot\|_L$ defined by
\begin{equation} \lb{2.12}
\|v_\bddot\|_L = \biggl(\, \sum_{n=1}^L \, \abs{v_j}^2\biggr)^{1/2}
\end{equation}
and $\langle \dott,\dott\rangle_L$ the associated inner product, we have that
\begin{align}
\|A_L(E)\|_{\HS}^2 &=
\|u(E,0)\|_L^2 \| u(E,\pi/2/)\|_L^2 - |\langle u(E,0), u(E,\pi/2)\rangle_L |^2
\lb{2.13} \\
&= \max_\theta \, \|u(E,\theta)\|_L^2 \, \min_\theta \|u(E,\theta)\|_L^2 \lb{2.14}
\end{align}
In particular,
\begin{equation} \lb{2.15}
\|A_L(E)\| \leq \sum_{j=0}^{L-1} \, \|T(j,E)\|^2
\end{equation}
\end{theorem}

\begin{remarks} 1. That \eqref{2.10}/\eqref{2.11} hold is either a direct calculation
verifying the formula or a calculation obtained by expanding $\binom{w_{n+1}(E')}
{w_n (E')}$ in terms of $\binom{u_{n+1}(E,\theta)}{u_n (E,\theta)}$ for $\theta =0,\pi/2$.
\eqref{2.13} is a direct calculation from \eqref{2.7} and \eqref{2.14} is a clever observation
\cite{KKL}.

2. Clearly, $\|T(n,E)\binom{1}{0}\|^2 = \abs{u_{n+1}(E)}^2 +
\abs{u_n(E)}^2 \geq \abs{u_n(E)}^2$, so $\|u(E,0)\|_L^2 \leq$ RHS of
\eqref{2.15}. Similarly for $\|u(E, \pi/2)\|$, so
$\|u(E,0)\|_L \|u(E, \pi/2)\|\leq$ RHS of \eqref{2.15}, so
\eqref{2.13} implies \eqref{2.15}.

3. $\|T(n)\|$ measures the growth of the fastest growing solution, so the RHS of \eqref{2.15}
in fact measures $(\max_\theta \|u(E,\theta)\|_L^2)$ and thus, by \eqref{2.14}, one could place
$\min_\theta \|u(E,\theta)\|/\max_\theta \|u(E,\theta)\|$ in front of the RHS of \eqref{2.15}.
\end{remarks}

Here is the key lower bound:

\begin{theorem}\lb{T2.2} Let $E',E''$ be two distinct zeros of $p_L(x)$. Then
\begin{equation} \lb{2.16}
\abs{E'-E_0} + \abs{E''-E_0} \geq \|A_L(E_0)\|^{-1}
\end{equation}
In particular,
\begin{equation} \lb{2.17}
\abs{z_L^{(+1)}(E_0) - z_L^{(-1)}(E_0)} \geq \|A_L (E_0)\|^{-1}
\end{equation}
\end{theorem}

\begin{proof} By \eqref{2.11}, for $E_j=E',E''$,
\begin{equation} \lb{2.18}
p_{\bddot -1}(E_j) =p_{\bddot -1}(E_0) + (E_j-E_0) A_L(E_0) p_{\bddot -1}(E_j)
\end{equation}
To say $p_L(E')=p_L(E'')=0$ says $E,E''$ are eigenvalues of
$J_{L;F}$ and $p_{\bddot -1}(\cdot)$
are the eigenvectors. So, by orthogonality of eigenvectors,
\begin{equation} \lb{2.19}
\langle p_{\bddot -1}(E'), p_{\bddot-1}(E'')\rangle_L=0
\end{equation}

By interchanging $E'$ and $E''$, if necessary, suppose
\begin{equation} \lb{2.20}
\|p_{\bddot-1}(E')\|_L\geq \|p_{\bddot -1}(E'')\|_L
\end{equation}
Take \eqref{2.18} for $E'$ and $E''$ and take the inner product with $p_{\bddot -1}(E')$
and subtract to get
\begin{align*}
\|p_{\bddot -1}(E')\|^2 &\leq \abs{E'-E_0}\, \abs{\langle p_{\bddot-1}(E'), A_L(E_0)
p_{\bddot -1}(E')\rangle} \\
& \qquad \qquad + \abs{E''-E_0}\, \abs{\langle p_{\bddot -1}(E'), A_L (E_0)
p_{\bddot -1}(E'')\rangle} \\
&\leq (\abs{E'-E_0} + \abs{E''-E_0}) \|A_L(E_0)\|\, \|p_{\bddot -1} (E')\|^2
\end{align*}
on account of \eqref{2.20}. \eqref{2.16} is immediate.

\eqref{2.13} follows from \eqref{2.16} and $z^{(-1)} \leq E_0 \leq z^{(+1)}$.
\end{proof}

While our main applications are to clock theorems and Poisson statistics, there is a
universal tunnelling bound.

\begin{theorem} \lb{T2.3} Let $J_{L;F}$ be a finite Jacobi matrix with $\alpha_-=
\inf a_n$, $\alpha_+=\max a_n$, $\beta=\max b_n -\min b_n$. Let
\begin{equation} \lb{2.21}
\gamma = \alpha_-^{-1} \abs{(\beta +2\alpha_+)^2 + \alpha_+^2 + 1}^{1/2}
\end{equation}
Then any pair of eigenvalues, $E,E'$, of $J_{L;F}$ obeys
\begin{equation} \lb{2.22}
\abs{E-E'} \geq \f{\gamma^2 -1}{\gamma^{2L} -1}
\end{equation}
\end{theorem}

\begin{remark} This bound is exponential, $\sim\gamma^{-2L}$, for $L$ large.
\end{remark}

\begin{proof} Adding a constant to $b_n$ does not change eigenvalue differences, so
we can suppose that
\begin{equation} \lb{2.23}
\max b_n =-\min b_n = \f{\beta}{2}
\end{equation}
Then any $E_0$ in the convex hull of $\spec (J_{L;F})$ obeys $\abs{E_0}\leq
\f{\beta}{2} + 2\alpha_+$, so $\abs{E_0-b_n}\leq \beta + 2\alpha_+$. Thus,
$\gamma$ is an upper bound on the Hilbert-Schmidt norm of
\[
\f{1}{a_{n+1}} \begin{pmatrix}
E-b_{n+1} & -a_n \\
1 & 0
\end{pmatrix}
\]
and so on the norm.

Pick $E$ between $E'$ and $E''$. It follows that
\[
\|T(j,E)\|\leq \gamma^j
\]
so \eqref{2.21} follows from \eqref{2.17} and \eqref{2.15}.
\end{proof}

One can also use our proof to see that one cannot have too many zeros near $E_0$.

\begin{theorem} \lb{T2.4} Define $\delta_{L,n}(E_0)$ to be the smallest number
so that
\begin{equation} \lb{2.24}
\#\{z_j^{(n)}\mid \abs{z_j^{(n)} -E_0} <\delta\}\leq n-1
\end{equation}
Then, for $n\geq 2$,
\begin{equation} \lb{2.25}
\delta_{L,n}(E_0) \geq \tfrac12\, \|A_L (E_0)\|_{\HS}^{-1} \, \sqrt{n}
\end{equation}
\end{theorem}

\begin{remarks} 1. If one has strong clock behavior, $\delta_{L,n}\sim cn/L$ for
$n$ fixed and $L$ large, so $\sqrt{n}$ is worse than one expects in nice cases.

2. Our proof shows that \eqref{2.17} can be ``improved," if $\|\cdot\|_{\HS}$
is used, to
\begin{equation} \lb{2.26}
\abs{z_L^{(+1)}(E_0) -z_L^{(-1)}(E_0)} \geq \f{1}{\sqrt{2}}\, \|A_L(E_0)\|_{\HS}^{-1}
\end{equation}
\end{remarks}

\begin{proof} There are at least $n$ zeros, $z_1, \dots, z_n$, in $\{z\mid \abs{z-E_0}
\leq\delta_{L,n}\}$. Order them so that if $\varphi_{j,m}=p_m (z_j)$, then $\|\varphi_j\|_L
\geq \|\varphi_2\|_L \geq \cdots \geq \|\varphi_n\|_L$. Let $\ti\varphi_j = \varphi_j/
\|\varphi_j\|_L$. Then the argument in the proof of Theorem~\ref{T2.2} says that,
for $j\leq k$,
\begin{align}
1 &\leq [\delta_{L,n} (\abs{A_{jj}} + \abs{A_{jk}})]^2 \notag \\
&\leq 2\delta_{L,n}^2 (\abs{A_{jj}}^2 + \abs{A_{jk}}^2) \lb{2.26x}
\end{align}
where
\begin{equation} \lb{2.27}
A_{jk} = \langle\varphi_j, A_L (E_0) \varphi_k\rangle
\end{equation}
and \eqref{2.26x} comes from $(a+b)^2 \leq 2(a^2 + b^2)$. Summing over all pairs and
noting that $\abs{A_{11}}^2$ occurs $n-1$ times, we find that
\[
\f{n(n-1)}{2} \geq 2\delta_{L,n}^2 (n-1) \|A\|_{\HS}^2
\]
which is \eqref{2.25}.
\end{proof}

\section{Clock Theorems for Bounded Variation Perturbations of
Free and Periodic OPRL} \lb{s3}

The basic result from which we will derive all our clock theorems presupposes
the existence of a complex solution to the difference equation \eqref{2.1} for
which we have precise information on the phase. The model is the Jost solution
which is asymptotic to $e^{i\theta (E)n}$, where $E=2\cos (\theta)$, $0\leq\theta
\leq \pi$.

\begin{theorem} \lb{T3.1} Let $\{a_n,b_n\}_{n=1}^\infty$ be a set of Jacobi parameters
and $\Delta$ a closed interval in $\bbR$. Suppose there exists a solution $u_n(E)$ of
\eqref{2.1} for $E\in\Delta$ which obeys
\begin{SL}
\item[{\rm{(i)}}]
\begin{equation} \lb{3.1}
 u_0(E)>0 \qquad \Ima u_1 (E)>0
\end{equation}
\item[{\rm{(ii)}}]
\begin{equation} \lb{3.1a}
u_n(E)=\abs{u_n(E)} \exp (i[n\alpha (E) + \beta_n(E)])
\end{equation}
where $\alpha$ is $C^1$ with
\begin{equation} \lb{3.2}
-\f{d\alpha}{dE} > 0 \qquad \text{all } E\in\Delta
\end{equation}
each $\beta_n$ continuous on $\Delta$, and
\begin{equation} \lb{3.3}
\lim_{n\to\infty} \, \biggl[\, \sup_{\substack{ E',E\in\Delta \\
\abs{E-E'} \leq 1/n}}  \abs{\beta_n(E') - \beta_n(E)}\biggr] =0
\end{equation}
\item[{\rm{(iii)}}] The transfer matrix $T_n(E)$ of \eqref{2.5} obeys
\begin{equation} \lb{3.4}
\tau \equiv \sup_{n, E\in\Delta}  \|T_n(E)\|<\infty
\end{equation}
\end{SL}
Then the density of states exists on $\Delta$,
\begin{equation} \lb{3.5}
d\nu(E)=-\f{1}{\pi}\, \f{d\alpha}{dE}\,dE
\end{equation}
and there is clock behavior uniformly in $\Delta$.
\end{theorem}

\begin{remarks} 1. $\alpha$ is, of course, a rotation number and \eqref{3.5} an
expression of the connection between the density of states and rotation numbers;
see Johnson-Moser \cite{JM82}.

2. \eqref{3.3} implies $\beta_n$ is irrelevant for eigenvalue spacing comparable
to $1/n$. To control possible spacings with $\Delta E$ small compared to $1/n$, one
needs some Lipschitz control of $\beta_n$, that is,
\[
\sup_{E\neq E'}\, \f{\abs{\beta_n(E') - \beta_n(E)}}{n\abs{E-E'}} \to 0
\]
which is where differentiability of $\beta_n$ and so moment conditions on $\{a_n,b_n\}$
came into \cite{Saff1}. We avoid this by using \eqref{3.4} to get a priori bounds.

3. \eqref{3.3} implies the same if $1/n$ is replaced by $A/n$ for any fixed $A$.
Define $\zeta_n (A)$ by
\begin{equation} \lb{3.5a}
\zeta_n(A) =\sup_{\substack{ E,E'\in\Delta \\ \abs{E-E'}\leq A/n}}
\abs{\beta_n(E')-\beta_n(E)}
\end{equation}

4. \eqref{3.3} is implied by an equicontinuity assumption, for example, uniform convergence
of $\beta_n$ to a continuous limit.
\end{remarks}

\begin{proof} By \eqref{3.1}, $u_n$ and $\bar u_n$ are independent solutions of \eqref{2.1}
and so cannot vanish at any points. Moreover,
\begin{equation} \lb{3.6}
p_{n-1} (E) = A(E)u_n(E) + \ol{A(E)}\,\, \ol{u_n(E)}
\end{equation}
where
\begin{equation} \lb{3.7}
A(E)= \f{W(p_{\bddot -1}, \bar u)}{W(u,\bar u)}
\end{equation}

Here, $W$ is the Wronskian. Given sequences, $f_n,g_n$,
\begin{equation} \lb{3.8}
W_n (f,g) = a_n (f_{n+1} g_n -f_n g_{n+1})
\end{equation}
is constant (call it $W(f,g)$) if $f,g$ both solve \eqref{2.1}. Since $p_{-1}=0$, $p_0=1$,
$a_0=1$, we have $W(p_{\bddot -1}, \bar u)=u_0(E)$, and clearly, $W(u,\bar u)=2 i\Ima (u_1
\bar u_0)=2i u_0 \Ima u_1$, so
\begin{equation} \lb{3.9}
A(E)=-\f{i}{2}\, \f{1}{\Ima u_1}
\end{equation}
is pure imaginary. Thus, \eqref{3.6} says $p_{n-1}(E)$ vanishes if and only if $u_n$ is real,
that is, by \eqref{3.1a},
\begin{equation} \lb{3.9b}
p_{n-1}(E)=0 \Leftrightarrow n\alpha (E) + \beta_n(E)=k\pi \qquad
(k\in\bbZ)
\end{equation}

Let
\begin{equation} \lb{3.9a}
\eta =\min_{E\in\Delta}\, \biggl[-\f{d\alpha}{dE}\biggr]
\end{equation}
Pick $N_0$ so that $n>N_0$ implies
\begin{equation} \lb{3.10}
\f{\eta}{\tau^2} \geq 4\zeta_n \biggl( \f{1}{\tau^2}\biggr)
\end{equation}
This can be done since $\zeta_n(A)\to 0$ as $n\to\infty$ by hypothesis. Since $\zeta_n(A)$
is increasing in $A$ and $\zeta_n(A+B)\leq\zeta_n(A) + \zeta_n(B)$, we have $\zeta_n (x\beta)
\leq \zeta_n([x]\beta) + \zeta_n (x\beta -[x]\beta) \leq ([x]+1)\zeta_n(\beta)\leq 2x
\zeta_n(\beta)$, if $x\geq 1$, \eqref{3.10} implies
\begin{equation} \lb{3.11}
q\geq \f{1}{\tau^2}\Rightarrow \eta q \geq 2\zeta_n(q)
\end{equation}

This in turn implies
\begin{equation} \lb{3.12a}
E -E'\geq \f{1}{n\tau^2} \Rightarrow \bigl[n[\alpha (E')-\alpha(E)]
+ [\beta_n(E') -\beta_n(E)]\bigr] \geq \f{n}{2} \, \eta\abs{E'-E}
\end{equation}

By \eqref{2.17} and \eqref{2.15}, any two successive zeros obey
\begin{equation} \lb{3.12}
\abs{E'-E}\geq \f{1}{n\tau^2}
\end{equation}
Thus, \eqref{3.12} implies that for $n > N_0$, any two solutions of \eqref{3.9b}
have distinct values of $k$. We also see from \eqref{3.12a} and continuity that if
$E$ is a solution of \eqref{3.9b}, there is another solution in $(E-\f{2\pi}{n\eta},E)$
and it has the next larger value of $k$ (i.e., $k+1$).

Subtracting \eqref{3.9b} for two successive values of \eqref{3.9b} and using $\zeta_n
(2\pi/ \eta)\to 0$, we see that
\begin{equation} \lb{3.13}
\sup_{\substack {E'<E \text{ successive} \\ \text{eigenvalues in }\Delta}}\,
\abs{n(\alpha(E') -\alpha(E)) -\pi}\to 0
\end{equation}
Given the uniformity of convergence of the difference quotient to the derivative,
\eqref{3.13} implies that
\begin{equation} \lb{3.14}
\sup_{\substack {E'<E \text{ successive} \\ \text{eigenvalues in }\Delta}}\,
\biggl| n(E'-E) \, \f{d\alpha}{dE}-\pi\biggr| \to 0
\end{equation}

This implies the density of states exists and is given by \eqref{3.5} and that one
has uniform clock behavior.
\end{proof}

\begin{theorem}[$\equiv$ Theorem~\ref{T1.3}] \lb{T3.2} Let $\{a_n,b_n\}_{n=0}^\infty$
be a set of Jacobi parameters obeying
\begin{equation} \lb{3.15}
\sum_{n=1}^\infty\, \abs{b_n} + \abs{a_n-1} <\infty
\end{equation}
so $\ess\,\supp(d\mu)=[-2,2]$. For any $\veps >0$, we have uniform clock behavior on
$[-2+\veps, 2-\veps]$.
\end{theorem}

\begin{remarks} 1. This includes Jacobi polynomials (rescaled to $[-2,2]$) for which
\[
\abs{b_n} + \abs{a_n-1}= O(n^{-2})
\]

2. Of course, the density of states is the free one.
\end{remarks}

\begin{proof} It is well-known (see, e.g., \cite{KS2,DSJost2}) that when \eqref{3.15} holds,
there exists, for all $x\in (-2,2)$, a solution $\ti u_n(x)$ so that if
\begin{equation} \lb{3.16}
z+z^{-1}=x
\end{equation}
(i.e., $z=e^{i\theta}$ with $x=2\cos\theta$ and $\theta\in (0,\pi)$), then
\begin{equation} \lb{3.17}
z^{-n} \ti u_n (x)\to 1
\end{equation}
uniformly on compact subsets of $(-2,2)$. Moreover, $\ti u_n(x)$ is continuous on $(-2,2)$
for each fixed $n$. By evaluating the Wronskian near $n=\infty$, we see
\begin{equation} \lb{3.18}
W(\ti u_n, \bar{\ti{u}}_n)=z-z^{-1}
\end{equation}
Thus, if
\begin{equation} \lb{3.19}
u_n (E)=\f{\ti u_n(E)}{\ti u_0(E)}
\end{equation}
then \eqref{3.1} and \eqref{3.1a} hold. If
\begin{equation} \lb{3.20}
u_0(E)=\abs{u_0(E)} e^{i\beta_\infty (E)}
\end{equation}
then we have that
\begin{equation} \lb{3.21}
\beta_n(E)\to\beta_\infty (E)
\end{equation}
uniformly on $[-2+\veps, 2-\veps]$. By continuity, \eqref{3.3} holds.

$\alpha$ is given by
\begin{equation} \lb{3.22}
\alpha(E) = \arccos \biggl(\f{E}{2}\biggr)
\end{equation}
so
\begin{equation} \lb{3.23}
-\f{1}{\pi}\, \f{d\alpha}{dE} = \f{1}{\pi}\, \f{1}{\sqrt{4-E^2}}
\end{equation}
and \eqref{3.2} holds.

Finally, standard variation of parameters about $z^n,z^{-n}$ shows that \eqref{3.4}
holds for each $\Delta =[2+\veps, 2-\veps]$.

Thus, Theorem~\ref{T3.1} applies, and we have clock behavior.
\end{proof}

In the above, we used the fact that $\beta_n\to\beta_\infty$ uniformly to obtain
\eqref{3.3}. In the bounded variation case, we will instead use:

\begin{lemma}\lb{L3.3} If $\beta_n =\beta_n^{(1)} + \beta_n^{(2)}$ where $\beta_n^{(1)}$
is $C^1$ and
\begin{equation} \lb{3.24}
\f{1}{n}\, \sup_{E\in\Delta}\, \biggl| \f{\partial\beta_n^{(1)}}{\partial E}\biggr|\to 0
\end{equation}
and $\beta_n^{(2)}\to\beta_\infty^{(2)}$ uniformly, then \eqref{3.3} holds.
\end{lemma}

\begin{proof} Immediate since $E,E'\in\Delta$ implies
\[
\abs{\beta_n (E)-\beta_n(E')} \leq \abs{E-E'}\, \sup_{E\in\Delta}\,
\biggl| \f{\partial\beta_n^{(1)}}{\partial x}\biggr| +
\abs{\beta_n^{(2)}(E) -\beta_n^{(2)}(E')}
\qedhere
\]
\end{proof}

\begin{theorem}[$\equiv$ Theorem~\ref{T1.4}]\lb{T3.4} Let $\{a_n,b_n\}_{n=0}^\infty$
be a set of Jacobi parameters with
\begin{equation} \lb{3.25}
a_n\to 1 \qquad b_n\to 0
\end{equation}
and
\begin{equation} \lb{3.26}
\sum_{n=1}^\infty (\abs{a_{n+1}-a_n} + \abs{b_{n+1}-b_n}) <\infty
\end{equation}
then for any $\veps >0$, we have uniform clock behavior in $[-2+\veps, 2+\veps]$.
\end{theorem}

\begin{remark} Again, the density of states is the free one by \eqref{3.25}.
\end{remark}

In order to prove this theorem, we need the following result:

\begin{theorem} \lb{T3.5} Let $B_n(\theta)$ depend continuously on $\theta\in I$,
a compact subinterval of $(0,\pi)$, and suppose
\begin{equation} \lb{3.27}
\sup_{\theta\in I}\, \sum_{n=1}^\infty\, \|B_{n+1}(\theta)-B_n(\theta)\| <\infty
\end{equation}
and
\begin{equation} \lb{3.28}
B_n(\theta)\to B_\infty(\theta)
\end{equation}
uniformly where $B_\infty(\theta)$ has eigenvalue $e^{\pm i\theta}$. Explicitly for
$V_\infty (\theta)$ continuous and invertible:
\begin{equation} \lb{3.28a}
B_\infty (\theta) = V_\infty(\theta) \begin{pmatrix} e^{i\theta} & 0 \\
0 & e^{-i\theta}\end{pmatrix} V_\infty (\theta)^{-1}
\end{equation}
Suppose each $B_n(\theta)$ has eigenvalues $e^{\pm i\varphi _n(\theta)}$ with
$\varphi_n (\theta)\in (0,\pi)$. Let
\begin{equation} \lb{3.29}
\ti T_n(\theta)=B_n(\theta)\dots B_1(\theta)
\end{equation}
Then
\begin{SL}
\item[{\rm{(i)}}]
\begin{equation} \lb{3.30}
\sup_{\theta\in I,n}\, \|\ti T_n(\theta)\| <\infty
\end{equation}
\item[{\rm{(ii)}}] There exists $S_\infty (\theta)$ depending continuously on $\theta$ so
that
\begin{equation} \lb{3.31}
\|\ti T_n(\theta) - V_\infty(\theta) D_n(\theta) S_\infty (\theta)\| \to 0
\end{equation}
where
\begin{equation} \lb{3.32}
D_n(\theta)=\begin{pmatrix} e^{i\eta_n(\theta)} & 0 \\
0 & e^{-i\eta_n (\theta)} \end{pmatrix}
\end{equation}
with
\begin{equation} \lb{3.33}
\eta_n(\theta)=\sum_{\eta=1}^n \varphi_n (\theta)
\end{equation}
\end{SL}
\end{theorem}

\begin{remark} See the notes to Section~2.1 of \cite{OPUC2} for a history of results on
bounded variation.
\end{remark}

\begin{proof} This is a strong version of Kooman's theorem \cite{Koo98}. In Section~12.1
of \cite{OPUC2}, \eqref{3.30} is proven, and in the notes to that section, it is noted that
\eqref{3.31} holds. In those notes, there is no $V_\infty$ because the analog of $D_n$
is not diagonal but can be diagonalized in a basis where $B_\infty$ is diagonal.
\end{proof}

\begin{proof}[Proof of Theorem~\ref{T3.4}] Define
\begin{equation} \lb{3.34}
B_n(\theta) = \begin{pmatrix}
(a_n a_{n+1})^{-1/2} (2\cos(\theta) -b_{n+1}) & -(a_n/a_{n+1})^{1/2} \\
(a_n/a_{n+1})^{-1/2} & 0 \end{pmatrix}
\end{equation}
so the transfer matrix at $E=2\cos\theta$ is
\begin{align}
T_n(\theta) &= \biggl( \f{a_n}{a_{n+1}}\biggr)^{1/2} B_n(\theta)
\biggl( \f{a_{n-1}}{a_n}\biggr)^{1/2} \dots B_1(\theta) \notag \\
&= a_{n+1}^{-1/2} \ti T_n(\theta) \lb{3.35}
\end{align}

Since $\det(B_n)=1$, $B_n(\theta)$ has eigenvalues $e^{\pm i\varphi_n}$,
$\varphi_n\in (0,\pi)$, if and only if
\begin{equation} \lb{3.36}
(a_n a_{n+1})^{-1/2} (2\cos\theta -b_{n+1}) = 2\cos(\varphi_n)\in (-2,2)
\end{equation}
and this holds uniformly for $\theta\in (\delta, \pi-\veps)$ and $n>N_0$ for
some fixed $N_0$. Thus, for $n\geq N_0$, we define
\begin{equation} \lb{3.37}
\eta_n(\theta) =\sum_{j=N_0}^n \varphi_j (\theta)
\end{equation}
and define
\begin{equation} \lb{3.38}
\ti D_n (\theta) = \begin{pmatrix}
e^{i\eta_n(\theta)} & 0 \\
0 & e^{-i\eta_n(\theta)} \end{pmatrix}
\end{equation}
and
\begin{equation} \lb{3.39}
\ti S_\infty (\theta) = S_\infty (\theta; B_n\dots B_{N_0}) B_{N_0-1} \dots B_1
\end{equation}
where $S_\infty (\theta;B_n \dots B_{N_0})$ is the $S_\infty$ for the sequence
$B_{N_0}, B_{N_0+1}, \dots$. Thus, \eqref{3.31} and \eqref{3.35} show that
\begin{equation} \lb{3.40}
\|T_n(\theta) -a_{n+1}^{-1/2} V_\infty (\theta) \ti D_n(\theta) \ti S_\infty (\theta)\|
\to 0
\end{equation}

We now proceed to construct a solution $u$ obeying the hypothesis
of Theorem~\ref{T3.1}.
Pick a two-vector $x(\theta)$ by $x(\theta)=e^{i\varphi_0(\theta)} \ti S_\infty
(\theta)^{-1} \binom{1}{0}$ where $\varphi_0$ is chosen below
and $u_n$ by
\begin{equation} \lb{3.41}
T_n (\theta) x(\theta) = \binom{u_{n+1}(\theta)}{u_n(\theta)}
\end{equation}
then \eqref{3.40} says that
\begin{equation} \lb{3.42}
u_n (\theta) = \abs{u_n(\theta)} \exp (i [\eta_n (\theta) + \beta_n^{(2)}(\theta)])
\end{equation}
where
\begin{equation} \lb{3.43}
\beta_n^{(2)}(\theta) \to \beta_\infty^{(2)}(\theta)
\end{equation}
uniformly.

Here $\beta_\infty^{(2)}(\theta)$ is $\varphi_0(\theta)$ plus the phase of
the $21$ element of $V_\infty(\theta)$, and so
$i(\eta_n(\theta) + \beta_\infty^{(2)}(\theta))$
is the phase of the lower component of
$a_{n+1}^{-1/2} V_\infty (\theta) \ti D_n (\theta) \ti S_\infty(\theta)x(\theta)$.
Since $u_n$ is not real, $u_0(\theta)\ne 0$, and so $\varphi_0(\theta)$ can be
chosen so that $u_0(\theta)> 0$.

Since $\varphi_j(\theta)\to\theta$, $\eta_{n+1}-\eta_n \to\theta$, and thus
the imaginary part of the Wronskian of $u$ and $\bar u$ is positive, so $u$ obeys
\eqref{3.1}.

By \eqref{3.36},
\begin{equation} \lb{3.44}
\f{\partial\varphi_n}{\partial E} =\f{\sin(\theta)}{\sin(\varphi_n)}\,
(a_n a_{n+1})^{-1/2}\, \f{\partial\theta}{\partial E}
\end{equation}
so
\begin{equation} \lb{3.45}
\f{\partial\eta_n}{\partial E} = \f{\partial\theta}{\partial E}\, (n-N_0) +
\sum_{j=N_0}^n \biggl( \f{\sin(\theta)}{\sin (\varphi_n)}\, (a_n a_{n+1})^{-1/2}
-1\biggr)
\end{equation}
Thus, if
\begin{equation} \lb{3.46}
\beta_n^{(1)} =\eta_n -n\, \f{\partial\theta}{\partial E}
\end{equation}
then
\begin{equation} \lb{3.47}
\f{1}{n}\, \f{\partial\beta^{(1)}}{\partial E} = \biggl( -\f{N_0}{n} + \f{1}{n}\,
\sum_{j=N_0}^n \biggl( \f{\sin(\theta)}{\sin(\varphi_n)}\, (a_n a_{n+1})^{-1/2}
-1\biggr)\biggr)\f{\partial\theta}{\partial E}
\end{equation}
converges uniformly to zero, since
\begin{equation} \lb{3.48}
\f{\sin(\theta)}{\sin(\varphi_n)}\to 1 \qquad (a_n a_{n+1})^{-1/2} \to 1
\end{equation}
uniformly in $\theta$.

Lemma~\ref{L3.3} applies, so condition (ii) of Theorem~\ref{T3.1} holds with
$\alpha(E)=\theta (E)=\arccos (\f{E}{2})$. We thus have clock behavior with density
of states the free one, that is, given by \eqref{3.23}.
\end{proof}

Now consider the periodic case, that is,
\begin{equation} \lb{3.49}
a_{n+p}=a_n \qquad b_{n+p}=b_p
\end{equation}
The spectrum now has bands (see, e.g., \cite{Last92}). For any $E$ in
the interior of the bands, there is a Floquet solution with
$u_{n+p}=e^{i\gamma(E)}u_n$ with $\gamma(E)\in (0,\pi)$ and
$\f{\partial\gamma}{\partial E}<0$. Thus
\begin{equation} \lb{3.50}
u_{Lp+r} = \abs{u_r} e^{i(L\gamma (E) +\beta_r^{(\infty)})}
\end{equation}
where $\beta_r^{(\infty)}$ is the phase of $u_r$. Theorem~\ref{T3.1} applies with
$\alpha =\gamma (E)/p$ and
\begin{equation} \lb{3.51}
\beta_{Lp+r}(E) = \beta_r^{(\infty)}(E) - \f{r}{p}\, \gamma (E)
\end{equation}
There are only $r$ such functions so \eqref{3.3} holds, and we recover the zero
spacing part of Theorem~2.6 of \cite{Saff3}.

If $a_n^{(0)}, b_n^{(0)}$ are periodic and $a_n =a_n^{(0)}+\delta a_n$, $b_n =
b_n^{(0)} +\delta b_n$ and
\begin{equation} \lb{3.52}
\sum_{n=1}^\infty (\abs{\delta b_n}+ \abs{\delta a_n}) <\infty
\end{equation}
then one can construct Jost solutions on the interiors of the bands. All that changes
is that \eqref{3.51} is replaced by
\begin{equation} \lb{3.53}
\lim_{L\to\infty}\, \beta_{Lp+r}(E) =\beta_r^{(\infty)}(E) - \f{r}{p}\, \gamma (E)
\end{equation}
so Theorem~\ref{T3.1} still applies. Similarly applying the ideas in the proof of
Theorem~\ref{T3.4}, we obtain a bounded variation result. Since it includes the
\eqref{3.52} result, we summarize in a single theorem:

\begin{theorem}\lb{T3.6} Let $a_n^{(0)}, b_n^{(0)}$ obey
\begin{equation} \lb{3.54}
a_{n+p}^{(0)}=a_n^{(0)}  \qquad b_{n+p}^{(0)} = b_n^{(0)}
\end{equation}
for some $p$. Let $a_n,b_n$ obey
\begin{gather}
\lim_{n\to\infty}\, \abs{a_n-a_n^{(0)}} + \abs{b_n-b_n^{(0)}} =0 \lb{3.55} \\
\sum_{n=1}^\infty \, (\abs{a_{n+p}-a_n} + \abs{b_{n+p} -b_n}) <\infty \lb{3.56}
\end{gather}
Then, for any closed interval which is a subset of the interior of the bands
{\rm{(}}see Remark{\rm{)}}, we have uniform clock behavior with density of states
$-\f{1}{p}\,\f{d\gamma}{dE}$.
\end{theorem}

\begin{remark} There are $p$ closed bands, $B_1,\dots, B_p$, generically disjoint
but perhaps touching in a single point (closed gap). By the interior of the bands,
we mean $\cup_{j=1}^p B_j^{\intt}$ which will be smaller than $(\cup_{j=1}^p
B_j)^{\intt}$ if some gap is closed, that is, we must remove all of the gaps,
including those that degenerate to single points.
\end{remark}

$p$ need not be the minimal period, so we have that

\begin{corollary}\lb{C3.7} Suppose
\[
a_n\to 1 \qquad b_n\to 0
\]
and, for some $p$, \eqref{3.56} holds. Then, in any closed interval in $\{E=2\cos
(\theta)\mid p\theta\neq 0,\, \f{\pi}{p}, \dots, \f{(p-1)\pi}{p}\}$, we have
uniform clock behavior.
\end{corollary}

As a final topic, we want to discuss zeros very near $E=2$ when
\begin{equation} \lb{3.57}
\sum_{n=1}^\infty n(\abs{a_n-1} + \abs{b_n})<\infty
\end{equation}
It should be possible to extend this argument to get uniform clock behavior in
$[-2,2]$, with a suitable modification to take into account the behavior exactly
at $\pm 2$. When \eqref{3.57} holds, the Jost function, $u$, can de defined on
$[-2,2]$; see, for example, the appendix to \cite{DSJost2}. If $u(2)=0$, we say
there is a resonance at $2$, and if $u(2)\neq 0$, we say that $2$ is nonresonant.

\begin{theorem} \lb{T3.8} Let $\{a_n, b_n\}_{n=1}^\infty$ be a set of Jacobi
parameters obeying \eqref{3.57}. Define $0\leq\theta_1^{(n)} < \theta_2^{(n)}
<\cdots$ so $E_j^{(n)}=2\cos(\theta_j^{(n)})$ are the zeros of $p_n (x)$ nearest
to $x=2$ and below. Then
\begin{SL}
\item[{\rm{(a)}}] If $2$ is a resonance, then
\begin{equation} \lb{3.58}
n\theta_j^{(n)}\to (j-\tfrac12)\pi
\end{equation}
\item[{\rm{(b)}}] If $2$ is nonresonant, then
\begin{equation} \lb{3.59}
n\theta_j^{(n)}\to j\pi
\end{equation}
\end{SL}
\end{theorem}

\begin{remarks} 1. The two simplest examples are the nonresonant $a_n\equiv 1$,
$b_n\equiv 0$ where
\begin{equation} \lb{3.60}
p_n (2\cos\theta) = c_n\, \f{\sin ((n+1)\theta)}{\sin\theta}
\end{equation}
and the resonant $a_n\equiv 1$ $(n\geq 2)$, $a_1=\sqrt2$, and $b_n\equiv 0$, where
\begin{equation} \lb{3.61}
p_n (2\cos\theta) = d_n \cos(n\theta)
\end{equation}

2. Notice that (for simplicity, consider the nonresonant case)
\[
E_j^{(n)} -E_{j+1}^{(n)} \sim \f{(2j+1)\pi^2}{n^2}
\]
the eigenvalue space is $O(n^{-2})$ and not equal in $E$ but has clock spacing in
$\theta$.
\end{remarks}

\begin{remark} The key fact that at a zero energy resonance, the scattering phase is
$\pi/2$ (mod $\pi$) and otherwise it is $0$ (mod $\pi$) is well-known in the
continuum case, for which there is extensive physics literature; see, for example,
Newton \cite{NewtonBk}.
\end{remark}

\begin{proof} By the theorems found in the appendix to \cite{DSJost2} (which
codifies well-known results), when \eqref{3.57} holds, one has (\cite[eqn.~(A.27)]{DSJost2})
\begin{equation} \lb{3.62}
\abs{p_n (e^{i\theta})}\leq C(n+1)
\end{equation}
and the existence of a solution $u_n (e^{i\theta})$ with
\begin{equation} \lb{3.63}
e^{-in\theta} u_n (e^{i\theta})\to 1
\end{equation}
uniformly on $\partial\bbD$. $u_0$ is called the Jost function and
\begin{equation} \lb{3.64}
W(p_{\bddot -1},u) =u_0
\end{equation}

We want to use \eqref{2.14} where there is a collision of notation, so we let
$v(\theta,\varphi)$ be the solution at $E=2\cos\theta$ and boundary condition $\varphi$.
Then \eqref{2.14} becomes
\begin{equation} \lb{3.65}
\|A_L(2\cos\theta) \| \leq \max_\varphi\, \|v(\theta,\varphi)\|_L
\min_\varphi\, \|v(\theta,\varphi)\|_L
\end{equation}
If $u_0 (\theta=0)\neq 0$, we get one solution for $\theta$ small, $v(\theta,
\varphi_0(\theta))$ which is uniformly bounded in $\theta$ and $n$, and another
solution $v(\theta,\varphi_1=0)$ ($=p_{\bddot-1}$) bounded by $Cn$. It follows that
\begin{equation} \lb{3.66}
\min_\varphi\, \|v(\theta,\varphi)\|_L\leq CL^{1/2} \qquad
\max_\varphi\, \|v(\theta,\varphi)\|\leq CL^{3/2}
\end{equation}
If $u_0(\theta=0)=0$, we start at $n=1$ since $u_1 (\theta=0)\neq 0$ and construct
the bounded and linearly growing solution that way (in essence, the two solutions
in this case are $p_{\bddot-1}$ and $q$, where $q$ is the second kind polynomial),
so \eqref{3.66} still holds.

We conclude, using Theorem~\ref{T2.2}, that when \eqref{3.57} holds, then
\begin{equation} \lb{3.67}
\sup_{\substack{E,E'\text{ successive zeros of }p_n \\ E,E'\in [2-\veps, 2]}}
\abs{E-E'}\geq \f{c}{n^2}
\end{equation}

Define $\varphi(\theta)$ by
\begin{equation} \lb{3.68}
u_0 (e^{i\theta})=\abs{u_0 (e^{i\theta})}e^{i\varphi (\theta)}
\end{equation}
$\varphi$ can be defined by $\theta\neq 0$ since $u_0(\theta)$ is then nonzero.
We can pick $\theta$ continuous on $(0,\veps)$. We claim
\begin{equation} \lb{3.69}
\varphi(0)\equiv\lim_{\theta\downarrow 0}\, \varphi(\theta) =
\begin{cases}
0 & \text{mod $\pi$ if $2$ is nonresonant} \\
\f{\pi}{2} &\text{mod $\pi$ if $2$ is a resonance}
\end{cases}
\end{equation}
Postponing the proof of this for now, let us complete the proof of the theorem.

By \eqref{3.64} and $W(u,\bar u)=z-z^{-1}$ (since $u\sim z^n$), we see that
\begin{align}
p_n (2\cos\theta) &= \f{\ol{u_0(e^{i\theta})}\, u_n (e^{i\theta}) -
u_0 (e^{i\theta})\, \ol{u_n (e^{i\theta})}}{2i\sin(\theta)} \lb{3.70} \\
&= \f{\abs{u_0 (e^{i\theta})}\, \abs{u_n (e^{i\theta})}}{2\sin\theta}\,
\sin(n\theta +\beta_n(\theta)) \lb{3.71}
\end{align}
where
\begin{equation} \lb{3.71a}
\beta_n(\theta) \to -\varphi(\theta)
\end{equation}
as $n$ goes to infinity uniformly in a neighborhood of $\theta =0$. Thus, zeros of
$p_n (2\cos\theta)$ are given as solutions of
\begin{equation} \lb{3.72}
n\theta + \beta_n (\theta) =j\pi
\end{equation}

In the resonant case, since
\begin{equation} \lb{3.73}
\sup_{\abs{\theta}\leq \f{c}{n}}\, \biggl|\beta_n (\theta) -
\f{\pi}{2}\biggr|_{\text{mod }\pi} \to 0
\end{equation}
there is at  least one solution asymptotic with
\begin{equation} \lb{3.74}
n\theta_j^{(n)} \to (j-\tfrac12)\pi
\end{equation}

If there were multiple solutions for some $j$, we would have two zeros with
$0<\theta' <\theta$,
\[
n(\theta -\theta')\to 0 \qquad
\theta < \f{\pi}{n}
\]
for $n$ large ($\pi$ can be any number strictly larger than $\pi/2$),
\begin{align*}
\abs{2\cos\theta' -2\cos\theta} &\leq \abs{\theta'-\theta}\sin(\theta) \\
&= o\biggl(\f{1}{n}\biggr) O \biggl(\f{1}{n}\biggr)
\end{align*}
violating \eqref{3.67}. Thus, there are unique solutions and \eqref{3.58} holds.

In the nonresonant case, \eqref{3.72} holds, but instead
\begin{equation} \lb{3.75}
\sup_{\abs{\theta} \leq \f{c}{n}}\, \abs{\beta_n(\theta)}_{\text{mod }\pi} \to 0
\end{equation}
which proves existence of solutions with
\begin{equation} \lb{3.76}
n\theta_j^{(n)}\to j\pi
\end{equation}
for $j=1,2,\dots$. We must prove uniqueness for $j\geq 1$ and nonexistence for
$j=0$.

The uniqueness argument for $j\geq 1$ is the same as in the resonant case. To show
no solution with $j=0$, we suppose that $J$ has $m$ eigenvalues above $E=2$ (by
Bargmann's bound \cite{HS02}, the number is finite). Let $J(\lambda)$ be the Jacobi
matrix with
\begin{align*}
a_n(\lambda) &= a_n \\
b_n(\lambda) &= \begin{cases}
b_n+\lambda & n\leq m+1 \\
b_n & n > m+1
\end{cases}
\end{align*}
It is easy to see that as $\lambda\to\infty$, $J(\lambda)$ has at least $m+1$
eigenvalues. So pick $\lambda_0=\inf \{\lambda\mid J(\lambda)$ has $m+1$ eigenvalues
in $(2,\infty)\}$. Then $\lambda_0 >0$ and $J(\lambda_0)$ has a resonance at $2$.
By the analysis of the resonant case, $p_n^{(\lambda_0)}(x)$ has $m$ zeros in
$(2,\infty)$ and its $(m+1)$st zero asymptotic to $2-(\f12\, \f{\pi}{n})^2$, which
means $p_n(x)$, whose zeros are less than those of $p_n^{(\lambda_0)}$, cannot have
a zero asymptotic to $\theta^{(n)}\to 0$.

That proves the result subject to \eqref{3.69}. In the nonresonant case, $u_0
(e^{i\theta})$ is continuous and nonvanishing at $\theta=0$, and $u_0(1)$ is real,
$\varphi(0)\equiv 0$ mod $\pi$, so continuity proves the top half of \eqref{3.69}.

In the resonant case, we note that a Wronskian calculation (see
\cite[eqn.~(A.49)]{DSJost2}) shows that
\begin{equation} \lb{3.77}
\Ima (u_1(e^{i\theta}) u_0(e^{i\theta}))=\sin\theta
\end{equation}
Since $u_0(1)=0$, $u_1(1)\neq 0$, and $u_1(1)$ is real, so
\begin{equation} \lb{3.78}
\lim_{\theta\downarrow 0}\, \Ima \biggl( \f{u_0(e^{i\theta})}{\sin\theta}\biggr)\neq 0
\end{equation}

On the other hand, $u_0 (e^{-i\theta})=\ol{u_0(e^{i\theta})}$, so $\varphi(\theta) +
\varphi(-\theta)\equiv 0$ mod $\pi$, which means that any limit point of $\varphi
(\theta)$ is a multiple of $\pi/2$. This is only consistent with \eqref{3.78} if the
limit is congruent to $\pi/2$ mod $\pi$.
\end{proof}

\section{Lower Bounds in the Ergodic Case} \lb{s4}

Our main goal in this section is to prove the following:

\begin{theorem}[$\equiv$ Theorem~\ref{T1.5}]\lb{T4.1} Let $a_n^{(\omega)}, b_n^{(\omega)}$
be ergodic Jacobi parameters. Let $E_0$ be such that
\begin{SL}
\item[{\rm{(i)}}] The Lyapunov exponent $\gamma (E_0)=0$.
\item[{\rm{(ii)}}] The symmetric derivative of $\nu$ exists at $E_0$ and is finite and nonzero.
\end{SL}
Then there exists $C>0$ so that with probability $1$,
\begin{equation} \lb{4.1}
\liminf_{n\to\infty}\, n [z_n^{(1)}(E_0) -z_n^{(-1)}(E_0)] >C
\end{equation}
\end{theorem}

This is particularly interesting because of the connection to Poisson behavior, where:

\begin{definition} We say a probabilistic family of Jacobi matrices has Poisson behavior
at $E_0$ if and only if for some $\lambda$ (normally $\lambda =$ density of zeros) we
have that for any $\alpha_1 < \beta_1 \leq \alpha_2 < \beta_2 \leq \cdots < \beta_\ell$
and any $k_1, k_2, \dots, k_\ell\in \{0,1,2,\dots\}$,
\begin{equation} \lb{4.2}
\begin{split}
\text{Prob} \biggl(\# \biggl\{ z_n^{(j)}(E_0) & \in \biggl[ E_0 + \f{\alpha_m}{n}, E_0 +
\f{\beta_m}{n}\biggr]\biggr\}
= k_m \text{ for } m=1, \dots, \ell \biggr)  \\
&\qquad \to \prod_{m=1}^\ell \f{[\lambda (\beta_m -\alpha_m)]^{k_m}}{k_m!} \,
e^{-\lambda (\beta_m-\alpha_m)}
\end{split}
\end{equation}
\end{definition}

\begin{remark} Poisson behavior was proven in certain random Schr\"odinger operators
by Molchanov \cite{Mol} and for random Jacobi matrices by Minami \cite{Min}. See
Stoiciu \cite{StoiJAT,StoiDiss} and Davies-Simon \cite{DavSim} for related work on OPUC.
\end{remark}

\begin{corollary}\lb{C4.2} Let $a_n^{(\omega)}, b_n^{(\omega)}$ be ergodic Jacobi
parameters and $E_0\in\bbR$ so that the symmetric derivative of $\nu$ exists at $E_0$
and is finite. Suppose there is Poisson behavior at $E_0$. Then $\gamma (E_0) >0$.
\end{corollary}

\begin{remarks} 1. Basically, \eqref{4.1} is a rigid level repulsion inconsistent
with Poisson behavior.

2. Ergodicity is critical here. Killip-Stoiciu \cite{KStoi} have examples which are
not ergodic for which there is Poisson behavior with $\gamma =0$.
\end{remarks}

\begin{proof}[Proof of Corollary~\ref{C4.2}] Suppose first $\gamma (E_0)=0$ so
Theorem~\ref{T4.1} applies. For each $n$, let $f_n(\omega)$ be the characteristic function
of $\{\omega\mid\abs{z_n^{(1)}(E_0) - z_n^{(-1)} (E_0)}\leq \f12 Cn^{-1}\}$. By \eqref{4.1},
$f_n(\omega)\to 0$ for a.e.\ $\omega$, so
\begin{equation} \lb{4.3}
\int f_n(\omega)\, d\omega\to 0
\end{equation}
as $n\to\infty$.

Clearly, if there is one $z^{(j)}$ in $[E_0 -\f14 Cn^{-1}, E_0]$ and one in $[E_0,
E_0 + \f14 Cn^{-1}]$, then $f_n (\omega)=1$. Thus, by the assumption of Poisson behavior,
\[
\lim \int f_n(\omega)\, d\omega \geq \biggl( \f{\lambda C}{4}\, e^{-\lambda C/4}\biggr)^2
\]
which contradicts \eqref{4.3}.
\end{proof}

Our proof of Theorem~\ref{T4.1} will use the complex solutions
constructed by Deift-Simon \cite{S169} and the estimate
Theorem~\ref{T2.2}. It is thus important to be able to estimate $A$
in terms of any pair of solutions with Wronskian $1$.

\begin{lemma}\lb{L4.3} Let $u^{(1)}, u^{(2)}$ be any pair of solutions of \eqref{2.1}
with Wronskian $1$. Then the kernel $K$ of \eqref{2.7} has the form
\begin{equation} \lb{4.4}
K(n,m,E) = u_n^{(1)} u_m^{(2)} - u_n^{(2)} u_m^{(1)}
\end{equation}
In particular,
\begin{equation} \lb{4.5}
\|A_L (E)\|_{\HS}^2 \leq \|u^{(1)}\|_L^2 \|u^{(2)}\|_L^2
\end{equation}
\end{lemma}

\begin{proof} Noting that $u_n^{(j)} u_m^{(j)} - u_m^{(j)} u_n^{(j)} =0$, we see that $K$
is invariant under linear changes of the $u$'s of determinant $1$. This proves \eqref{4.4}.
\eqref{4.5} follows as in Theorem~\ref{T2.1}.
\end{proof}

We need the following result of Deift-Simon \cite{S169}:

\begin{theorem}\lb{T4.4} Let $a_n^{(\omega)},b_n^{(\omega)}$ be ergodic Jacobi parameters and let
{\rm{(i)--(ii)}} of Theorem~\ref{T4.1} hold for $E_0$. Then for a.e.\ $\omega$, there exists
a complex-valued solution $u(\dott, \omega)$ of \eqref{2.1} so that
\begin{SL}
\item[{\rm{(i)}}] The Wronskian of $u$ and $\bar u$ is $-2i$.
\item[{\rm{(ii)}}] We have that
\begin{equation} \lb{4.6}
\limsup_{n\to\infty}\, \abs{n}^{-1} \sum_{j=0}^{n-1} \,
\abs{u(j,w)}^2 \leq 2 \pi\, \f{d\nu}{dE}
\end{equation}
\end{SL}
\end{theorem}

\begin{remarks} 1. In fact, $\abs{u(j,w)} = \abs{u(0, T^j w)}$ and $\bbE (\abs{u(0,w)}^2)
<\infty$, so by the Birkhoff ergodic theorem, the $\limsup$ on the left of \eqref{4.6}
can be replaced by a limit which is a.e.\ constant with a constant bounded by the right
side of \eqref{4.6}

2. \cite{S169} states results for a.e.\ $E$ with $\gamma(E)=0$, but the proof shows
that what is needed is (i)--(ii).
\end{remarks}

\begin{proof}[Proof of Theorem~\ref{T4.1}] By \eqref{4.5}, taking $u^{(1)}=u$ and
$u^{(2)}=(2i)^{-1}\bar u$, we see
\begin{equation} \lb{4.7}
\limsup\, n^{-1} \|A_n(E)\| \leq \pi\, \f{d\nu}{dE}
\end{equation}
by \eqref{4.6}. Thus, by Theorem~\ref{T2.2},
\[
\liminf_{n\to\infty}\, n [z_n^{(1)}(E_0)- z_n^{(-1)}(E_0)]
\geq \f{1}{\pi}\left[\f{d\nu}{dE}\right]^{-1}
\]
which is \eqref{4.1}.
\end{proof}

\section{Upper Bounds via Transfer Matrices} \lb{s5}

Our goal in this section is to prove:

\begin{theorem}\lb{T5.1} Let $\{a_n,b_n\}_{n=1}^\infty$ be a set of Jacobi parameters.
For any bounded interval $I\subset\bbR$, we have
\begin{equation} \lb{5.1}
\sup_{\substack{E,E' \text{ successive} \\ \text{zeros of $p_n$ in $I$}}}
\abs{E-E'} \leq \f{[8e\sup_{E\in I} \|T_n(E)\|] (\prod_{j=1}^n a_j)^{1/n}}{n}
\end{equation}
\end{theorem}

\begin{corollary}\lb{C5.2} Let $\{a_n,b_n\}_{n=1}^\infty$ be a set of Jacobi parameters
and $I=[\alpha,\beta]$ a closed interval. Let
\begin{equation} \lb{5.2}
A=\sup_n a_n <\infty \qquad
T=\sup_{n;\,E\in I}\, \|T_n(E)\|<\infty
\end{equation}
Let $c=8eAT$. Then for any $\delta >0$, there is $N_0$ so that if $n>N_0$ and
$E\in [\alpha+\delta, \beta-\delta]$ is a zero of $p_n$, then there are at least
two additional zeros in $[E-\f{c}{n}, E+\f{c}{n}]$, one above $E$ and one below.
\end{corollary}

\begin{proof}[Proof of Corollary~5.2] It is known that any point in $\spec(J)$ is a
limit point of zeros, so for large enough $N_0$, there are zeros in $[\alpha,
\alpha+\delta)$ and $(\beta-\delta,\beta]$. Thus, \eqref{5.1} implies the result.
\end{proof}

\begin{example}\lb{E5.3} Let $a_n\equiv 1$, $b_n\equiv 0$ so
\begin{equation} \lb{5.3}
p_n (2\cos\theta) =\f{\sin((n+1)\theta)}{\sin\theta}
\end{equation}
and for $n$ odd, $p_n(0)=0$. The next nearest zero is at $\theta
=\f{\pi}{2} - \f{\pi}{n+1}$, so at $E\sim\f{2\pi}{n+1}$
($\left.\f{d}{d\theta}(2\cos\theta) \right|_{\theta=\f{\pi}{2}}
=-2$). In this case, $p_n$ at $E=0$ ($\theta=\f{\pi}{2}$) is
$(1,0,-1,0,1,0, \dots)$ for $n=0,1,2,\dots$ and $q_n =
(0,-1,0,1,\dots)$, so $T_{n\text{ odd}}=\left(\begin{smallmatrix} 0
& +1 \\ \mp 1 & 0 \end{smallmatrix} \right)$ so $\|T_n(0)\|=1$ and
$\|T_n(E)\|\sim 1$ for $E$ near $0$. Thus, the correct answer for
the spacing is $2\pi\sim 6.3$ and our upper bound is $8e\sim 21.7$,
a factor of about $3.5$ too large. \qed
\end{example}

To get Theorem~\ref{T5.1}, we will use

\begin{theorem}\lb{T5.4} Let $Q$ be a polynomial with all its zeros real. Let $Q(E_0)=0$,
$Q'(E_1)=0$, with $E_0<E_1$ and $Q$ nonvanishing on $(E_0,E_1)$. Then
\begin{equation} \lb{5.4}
\abs{E_1-E_0} \leq e\, \f{\abs{Q(E_1)}}{\abs{Q'(E_0)}}
\end{equation}
\end{theorem}

\begin{proof} Since $\f{d^2}{dx^2}\log (x-x_0) = -\f{1}{\abs{x-x_0}^2}<0$, we see that
$g(E)=\log (\abs{Q(E)}/\abs{E-E_0})$ is concave on $[E_0,E_1]$. Note that $g(E_0) =
\log\abs{Q'(E_0)}$ and that the tangent to $g$ at $E_1$ is
\begin{equation} \lb{5.5}
\log\biggl( \f{\abs{Q(E_1)}}{\abs{E_1-E_0}}\biggr) - \f{1}{\abs{E_1-E_0}}\, (E-E_1)
\end{equation}
Thus,
\[
\abs{Q'(E_0)}\leq \biggl( \f{\abs{Q(E_1)}}{\abs{E_1-E_0}}\biggr) e
\]
which is \eqref{5.4}.
\end{proof}

The $Q$ we will take to get \eqref{5.1} is not $P$, but the discriminant
\begin{equation} \lb{5.6}
\Delta_n (E)=\tr (T_n(E))
\end{equation}
associated to the periodic set of Jacobi parameters
\begin{equation}\lb{5.7}
a_{mn+q}^{(n)}=a_q \qquad b_{mn+q}^{(n)}=b_q \qquad q=1,\dots, n;\, m\geq 0
\end{equation}
We have:

\begin{lemma}\lb{L5.5} The zeros of $p_{n-1}$ and $\Delta_n$ interlace. Thus, if $E_1 <
E_2 < E_3$ are three successive zeros of $\Delta_n$, then $p_{n-1}$ has two zeros,
$E$ and $E'$, with
\begin{equation} \lb{5.8}
\abs{E-E'} < \abs{E_3 - E_1}
\end{equation}
\end{lemma}

\begin{proof} We need the analysis of $\Delta_n$ as a periodic discriminant \cite{Last92}.
$\Delta_n$ has $n$ bands given by $\alpha_1 < \beta_1 \leq \alpha_2 < \beta_2 \leq
\cdots \leq \alpha_n < \beta_n$ and bands $[\alpha_j,\beta_j]$. $\Delta^{-1}
(-2,2)=\cup_{j=1}^n (\alpha_j,\beta_j)$ and, in particular, the zeros of $\Delta$
lie one per band. $p_{n-1}$ has one zero in each gap $[\beta_1,\alpha_2], [\beta_2,
\alpha_3], \dots, [\beta_{n-1},\alpha_n]$. That gives us the interlacing. \eqref{5.8}
is an immediate consequence of this interlacing.
\end{proof}

To get a bound on $\Delta'_n$ at its zeros, we need a bound on the rotation number
for ergodic Schr\"odinger operators found by Deift-Simon \cite{S169}. This rotation
number is
\begin{equation} \lb{5.9}
\alpha(E)=\pi (1-\nu (-\infty,E))
\end{equation}
where $\nu$ is the density of states. Thus, $\alpha$ runs from $\pi$ to $0$ as $E$
runs from $\min\spec(J)$ to $\max\spec (J)$. $\cos(\alpha)$ runs from $-1$ to $1$.

\begin{proposition}\lb{P5.6} In the periodic case, on $\spec (J)$,
\begin{equation} \lb{5.10}
\f{d\cos(\alpha(E))}{dE} \geq \f12 \biggl(\, \prod_{j=1}^n a_j\biggr)^{-1/n}
\end{equation}
\end{proposition}

\begin{remarks} 1. (1.2) of \cite{S169} is an integrated form of \eqref{5.10}. We can
take derivatives since $\spec(J)$ is a union of intervals.

2. Deift-Simon assume $a_j\equiv 1$. By using the modification of the Thouless
formula for general $a_j$, it is easy to see their proof yields \eqref{5.10}.

3. In the free case ($a_n$ constant, $b_n=0$), one has equality in \eqref{5.10}.
\end{remarks}

\begin{proof}[Proof of Theorem~\ref{T5.1}] $\Delta_n(E)$ is related to $\alpha (E)$
by
\begin{equation} \lb{5.11}
\Delta_n (E)=2 \cos (n\alpha (E))
\end{equation}
Thus, at zeros of $\Delta_n$ where $\cos(n\alpha(E))=0$, we have $\abs{\sin (n\alpha(E))}
=1$. So at such points,
\begin{align*}
\biggl| \f{d\Delta_n}{dE}\biggr| &= 2n\, \f{d\alpha}{dE} \geq 2n \sin(\alpha(E))\, \f{d\alpha}{dE} \\
&= 2n\, \f{d\cos(\alpha(E))}{dE} \\
&\geq n \biggl(\, \prod_{j=1}^n a_j\biggr)^{-1/n}
\end{align*}
by \eqref{5.10}.

Since $\Delta_n =\tr (T_n)$, $\abs{\Delta_n}\leq 2\|T_n\|$, so \eqref{5.4} becomes
\begin{equation} \lb{5.12}
\abs{E_1-E_0} \leq 2e \, \f{\|T_n\| (\prod_{j=1}^n a_j)^{1/n}}{n}
\end{equation}

Between the zeros of $\Delta$ are two $(E_0,E_1)$-type intervals and
so, between the first and third of three zeros are four such
intervals. \eqref{5.8} and \eqref{5.12} imply \eqref{5.1}.
\end{proof}


\section{Pr\"ufer Angles and Bounds on Zero Spacing} \lb{s6}

There are various possible Pr\"ufer angles. We will exploit one that is ideal
for studying the energy dependence of zeros of $p_n$.

\begin{proposition}\lb{P6.1} Fix Jacobi parameters $\{a_n,b_n\}_{n=1}^\infty$. For
each $n=1,2,\dots$, there is a unique continuous function $\theta_n (E)$ determined by
\begin{align}
\tan (\theta_n(E)) &= \f{p_n(E)}{p_{n-1}(E)}  \lb{6.1} \\
\lim_{E\to -\infty}\, \theta_n(E) &= -\f{\pi}{2} \lb{6.2}
\end{align}
Moreover,
\begin{equation} \lb{6.3}
\f{d\theta_n}{dE} = \f{\sum_{j=0}^{n-1} p_j (E)^2}{a_n (p_{n-1} (E)^2 + p_n(E)^2)} >0
\end{equation}
\end{proposition}

\begin{proof} $p_n(E)/p_{n-1}(E)\to -\infty$ as $E\to -\infty$ and the ratio is
continuous on $\bbR\cup\{\infty\}$, so existence and uniqueness are immediate,
as is \eqref{6.2} since $\tan (-\pi/2)=-\infty$. Note next that since $\f{d}{dy}
\arctan (y)=\f{1}{1+y^2}$, we have
\begin{equation} \lb{6.4}
\f{d\theta_n}{dE} = \f{p_{n-1} p'_n-p_n p'_{n-1}}{p_{n-1}^2 (1+p_n^2/p_{n-1}^2)}
\end{equation}
so that \eqref{6.3} follows from the CD formula \eqref{A.1.3}.
\end{proof}

\begin{remark} \eqref{6.3} is closely related to well-known formulae relating the
derivative of projective angles in $\bbS\bbU (1,1)$ or $\bbS\bbL(2,\bbR)$;
see Theorem~10.4.11 of \cite{OPUC2}. These formulae have been used by Furstenberg,
Carmona, and others; see, for example, Carmona \cite{Car82} or Simon
\cite{Molchproc}.
\end{remark}

The transfer matrix, $T_n(E)$, is a $2\times 2$ matrix of determinant $1/a_{n+1}$, so
\begin{equation} \lb{6.5}
\|T_n(E)^{-1}\| = a_{n+1}\|T_n(E)\|
\end{equation}
Thus, since $(p_{j+1}, p_j)^t =T_j(E) (1,0)^t$, we have that
\begin{equation} \lb{6.6}
a_{j+1}^{-2}\|T_j(E)\|^{-2} \leq p_{j+1}(E)^2 + p_j(E)^2 \leq \|T_j(E)\|^2
\end{equation}

Since
\begin{equation} \lb{6.7}
\tfrac12 \sum_{j=0}^{n-2} \, (p_j^2 + p_{j+1}^2) \leq \sum_{j=0}^{n-1} p_j^2
\leq \sum_{j=0}^{n-2} \, (p_j^2 + p_{j+1}^2)
\end{equation}
\eqref{6.3} immediately implies

\begin{theorem}\lb{T6.2} Let
\begin{equation} \lb{6.8}
t_n(E) =\sup_{0\leq j\leq n-1} \, (1+a_{j+1}^2)\|T_j(E)\|^2
\end{equation}
Then
\begin{equation} \lb{6.9}
\frac{n}{2t_n(E) t_{n-1}(E)}
\leq a_n \, \f{d\theta_n}{dE}
\leq nt_n(E) t_{n-1}(E)
\end{equation}
\end{theorem}

This, in turn, implies $1/n$ upper and lower bounds on zero spacings sufficient
for what we needed in Section~\ref{s3}:

\begin{theorem}\lb{T6.3} If $\Delta$ is an interval in $\bbR$ on which $\tau_\Delta
\equiv \sup_{E\in\Delta, n}\, (1+a_{n+1}^2)\|T_n(E)\|^2<\infty$, then
\begin{equation} \lb{6.10}
\inf_{\substack{ E,E'\in\Delta \\ E,E'\, \text{successive zeros of $p_n(E)$}}}
\abs{E-E'} \geq \f{a_n \pi}{\tau_\Delta^2 n}
\end{equation}
and
\begin{equation} \lb{6.11}
\sup_{\substack{ E,E'\in\Delta \\ E,E'\, \text{successive zeros of
$p_n(E)$}}} \abs{E-E'} \leq \f{2a_n\pi \tau_\Delta^2}{n}
\end{equation}
Moreover, if $\Delta=[\alpha,\beta]$, $p_n$ has zeros in $[\beta -\f{2a_n\pi\tau_\Delta^2}
{n},\beta]$ and $[\alpha, \alpha + \f{2a_n\pi \tau_\Delta^2}{n}]$ once $\abs{\beta-\alpha}
\geq \f{2a_n \pi\tau_\Delta^2}{n}$.
\end{theorem}

\begin{proof} Since $\theta_n$ is monotone in $E$ and $p_n(E)=0$ if and only if
$\theta_n(E) =\ell\pi$ for some $\ell\in\bbZ$, we have at successive zeros, $E<E'$,
that
\begin{equation} \lb{6.12}
\pi = \theta_n(E') -\theta_n(E) =\int_E^{E'}\, \f{d\theta_n}{dE}\, dE
\end{equation}
\eqref{6.10}/\eqref{6.11} are then immediate from \eqref{6.9}. The final assertion
comes from the fact that $\theta (E_1)-\theta (E_0)\geq \pi$ implies that $\tan
(\theta(E))$ has a zero in $[E_0,E_1]$.
\end{proof}

\section{Relations of the Weight to the Christoffel Function} \lb{s7}

The previous sections were dominated by the transfer matrix. In this section, we shift
to the weight where the CD kernel (see \eqref{A.1.1}) will play a major role. This section
is a technical interlude: a detailed result that will be useful in the analysis of
examples in later sections. Our main result in this section is

\begin{theorem}\lb{T7.1} Suppose that
\begin{equation} \lb{7.1}
d\mu = w(x)\, dx + d\mu_\s
\end{equation}
with $d\mu_\s$ singular, and that for some $x_0$, $a>0$, and some $r>0$,
\begin{equation} \lb{7.2}
\int_{x_0-a}^{x_0+a} w(x)^{-r}\, dx < \infty
\end{equation}
Then
\begin{equation} \lb{7.3}
K_n(x_0,x_0) \leq C_r n^{1+r^{-1}}
\end{equation}
where $C_r$ only depends on $r,a$ and the integral in \eqref{7.2}.
\end{theorem}

This result generalizes one of Geronimus (see \cite[Remark~3.3 and Table~II]{GBk})
in two ways. His estimate is on $\abs{\varphi_n}^2$ not $K_n$ and, more importantly,
his estimates require global estimates on $w$ in the context of OPUC rather than
just our local estimate. One reason we can go beyond Geronimus is that he uses the
Szeg\H{o} function and we just use the Christoffel variational principle. Another
reason is that we have a powerful result of Nevai \cite{Nev79}:

\begin{proposition}[\cite{Nev79}]\lb{P7.2} For any $p$ in $(0,\infty)$, there is a
constant $D_p$ so
\begin{equation} \lb{7.4}
\int_{-1}^1 \abs{\pi_n(x)}^p\, dx \geq D_p n^{-1} \abs{\pi_n(0)}^p
\end{equation}
for any polynomial $\pi_n$ of degree $n$.
\end{proposition}

\begin{proof} Since this is a special case of Nevai's result which depends on
several arguments, for the reader's convenience, we extract exactly what is needed
for \eqref{7.4}.

Let
\begin{equation} \lb{7.5}
d\mu_0 (x) =\f{1}{\pi}\, (1-x^2)^{-1/2}\, dx
\end{equation}
on $[-1,1]$, so if $x=\cos(\theta)$, then
\[
d\mu_0=\f{d\theta}{\pi}
\]
on $[0,\pi]$ which implies, as is well known, that the OPs for \eqref{7.5} are given by
\begin{equation} \lb{7.6}
p_n (\cos(\theta)) = \begin{cases}
1 & n=0 \\
\sqrt2\, \cos (n\theta) & n\geq 1
\end{cases}
\end{equation}
the Chebyshev polynomial (of the first kind).

It follows that for $x=\cos(\theta)\in [-1,1]$,
\begin{equation} \lb{7.7}
K_n (x,x;d\mu_0)=1 + \sum_{j=1}^n 2\cos^2 (jn) \leq 2n+1
\end{equation}
Thus, by \eqref{A.1.7},
\begin{equation} \lb{7.8}
\sup_{\abs{x}\leq 1}\, \abs{\pi_n(x)}^2 \leq (2n+1) \int_{-1}^1 \abs{\pi_n(x)}^2
(1-x^2)^{-1/2} \, \f{dx}{\pi}
\end{equation}

If $m$ is an integer, $(\pi_n)^m$ is a polynomial of degree at most $mn$, so
\eqref{7.8} implies
\begin{equation} \lb{7.9}
\sup_{\abs{x}\leq 1}\, \abs{\pi_n(x)}^{2m} \leq (2mn+1) \int_{-1}^1 \abs{\pi_n(x)}^{2m}
(1-x^2)^{-1/2} \, \f{dx}{\pi}
\end{equation}
If $2m-2 < p\leq 2m$, we write
\[
\abs{\pi_n}^{2m} \leq \abs{\pi_n}^p \bigl(\, \sup_{\abs{x}\leq 1}\,
\abs{\pi_n(x)}\bigr)^{2m-p}
\]
to deduce
\begin{equation} \lb{7.10}
\sup_{\abs{x}\leq 1}\, \abs{\pi_n(x)}^p \leq \biggl\{2n \biggl(\biggl[ \f{p}{2}\biggr]
+1\biggr)\biggr\} \int_{-1}^1 \abs{\pi_n(x)}^p (1-x^2)^{-1/2} \,\f{dx}{\pi}
\end{equation}

Given $p$, pick $\ell$ from $1,2,\dots$ so $\ell p\geq\f12$ and apply \eqref{7.10} to the
polynomial $\pi_n(x)(1-x^2)^\ell$ which has degree $n+2\ell$ and get
\begin{equation} \lb{7.11}
\abs{\pi_n(0)}^p \leq \sup_{\abs{x}\leq 1}\, \abs{(1-x^2)^\ell
\pi_n(x)}^p \leq \biggl\{2 (n+2\ell) \biggl(
\biggl[\f{p}{2}\biggr] +1 \biggr)\biggr\} \int_{-1}^1
\abs{\pi_n(x)}^p \, \f{dx}{\pi}
\end{equation}
since $(1-x^2)^{\ell p-\f12} \leq 1$.

Find $D_p$ so for $n\geq 1$,
\[
2(n+2\ell) \biggl( \biggl[\f{p}{2}\biggr] +1\biggr) \leq D_p^{-1} n
\]
and \eqref{7.11} implies \eqref{7.4}.
\end{proof}

\begin{proof}[Proof of Theorem~\ref{T7.1}] By \eqref{A.1.15}, we can suppose $d\mu_\s
=0$. By scaling and translation, we can suppose $x_0=0$, $a=1$. By Theorem~\ref{T.A.1},
we need to get lower bounds on $\int_{-1}^1 \abs{\pi_n(x)}^2 w(x)\, dx$. By H\"older's
inequality, for any $\alpha,\beta$, $p\in (1,\infty)$, and $q$ dual to $p$,
\begin{align}
\int_{-1}^1 \abs{\pi_n(x)}^\alpha\, dx
&= \int_{-1}^1 \abs{\pi_n(x)}^\alpha w(x)^\beta w(x)^{-\beta}\, dx \notag \\
&\leq \biggl( \int_{-1}^1 \abs{\pi_n(x)}^{\alpha p} w(x)^{\beta p}\biggr)^{1/p}
\biggl( \int_{-1}^1 w(x)^{-\beta q} \, dx\biggr)^{1/q} \lb{7.12}
\end{align}
We want to pick $\beta,q,\alpha$ so $\beta q=r$, $\alpha p=2$, $\beta p=1$, that is,
\begin{equation} \lb{7.13}
q=1+r \,\qquad p=\f{1+r}{r} \,\qquad \alpha =\f{2r}{1+r} \,\qquad \beta =\f{r}{1+r}
\end{equation}

The result is that
\begin{align}
\int_{-1}^1 \abs{\pi_n(x)}^2 w(x)\, dx
&\geq \biggl( \int_{-1}^1 w(x)^{-r}\, dx\biggr)^{-1/r} \biggl( \int_{-1}^1
\abs{\pi_n(x)}^{2r/(1+r)}\, dx\biggr)^{(1+r)/r} \notag \\
&\geq C\abs{\pi_n(0)}^2 n^{-1-r^{-1}} \lb{7.14}
\end{align}

Taking the $\inf$ over all $\pi_n$'s with $\pi_n(0)=1$ and using \eqref{A.1.5}, we get
\[
K_n (0,0)^{-1} \geq Cn^{-1-r^{-1}}
\]
which is \eqref{7.3}.
\end{proof}

\begin{example}\lb{E7.3} Let $d\mu$ be the measure on $[-1,1]$ given by
\begin{equation} \lb{7.14a}
d\mu(x)= C_{a,b} \abs{x}^a (1-\abs{x}^2)^b\, dx
\end{equation}
where $a\geq 0$, $b\geq -1$, and $C$ is a normalization constant.

This is an even measure so $p_{2n-1}(0)=0$. Moreover,
\begin{equation} \lb{7.14b}
p_{2n}(x) = q_n (x^2)
\end{equation}
where $q_n$ are the OPs for the measure obtained from an $x\to y=x^2$ change of
variables. Since $dx=(dy)/y^{1/2}$, we see $q_n$ are the orthogonal polynomials
for the measure
\begin{equation} \lb{7.14c}
\ti C_{a,b} \abs{y}^{a/2-1/2} (1-y)^b \, dy
\end{equation}

Thus, up to a constant,
\[
K_n (0,0;d\mu) = \ti K_{n,\alpha.\beta}(1,1)
\]
where $\ti K_{n;\alpha,\beta}$ is the CD kernel for the Jacobi polynomial associated
to $(1-x)^\alpha (1+x)^\beta\, dx$ with $\alpha = \f12 a -\f12$, $\beta =b$. If we
call these orthogonal polynomials $j_{\alpha,\beta}$, and $J_{\alpha,\beta}$ the
conventional normalization, then \cite{Szb,MWorld}
\begin{align*}
\|J_{\alpha,\beta}\|_2 &\sim C_{\alpha,\beta}^{(1)} n^{-1} \\
J_{\alpha,\beta}(1) &\sim C_{\alpha,\beta}^{(2)} n^\alpha
\end{align*}
so
\[
j_{\alpha,\beta}(1) \sim (C_{\alpha,\beta}^{(1)})^{-1/2} (n^{1/2}) J_{\alpha,\beta}(1)
= C_{\alpha,\beta}^{(3)} n^{\alpha + 1/2}
\]
and
\[
\ti K_{n;\alpha,\beta}(1,1) \sim n^{2\alpha +2}
\]

Taking $\alpha = \f12 (a-1)$, we get
\begin{equation} \lb{7.14d}
K_n (0,0;d\mu) \sim n^{1+ a}
\end{equation}
We can take $r$ in Theorem~\ref{T7.1} arbitrary with $ra<1$, so \eqref{7.3}
cannot be improved.
\qed
\end{example}

The following shows that in some cases the power of $n$ in Theorem~\ref{T7.1} is optimal:

\begin{theorem}\lb{T7.4} Let $d\mu(x)=w(x)\, dx$ where $\supp (d\mu)\subset [-1,1]$ and
\begin{equation} \lb{7.25}
\abs{w(x)} \leq C\abs{x}^\alpha
\end{equation}
for some $\alpha <1$. Then
\begin{equation} \lb{7.26}
\abs{K_n(0,0)} \geq C_1 n^{1+\alpha}
\end{equation}
\end{theorem}

\begin{remark} For $w(x)=C_\alpha \abs{x}^\alpha$ on $[-1,1]$, \eqref{7.2} holds for any
$r<1/\alpha$, so \eqref{7.26} says \eqref{7.3} cannot hold for any smaller power of $n$
in case $r>1$.
\end{remark}

\begin{proof} Let $\pi_n$ be the polynomial of Theorem~\ref{T.A.6} where $x_0=0$,
$a=1$. On $\abs{x}\in [\f{j}{n}, \f{j+1}{n}]$, $j=0,1,2,\dots, n-1$,
\[
\abs{\pi_n(x)} \leq \begin{cases}
1 & j=0 \\
\f{1}{2n} + \f{1}{2j} & j=1,2,\dots
\end{cases}
\]
so, by \eqref{7.25},
\begin{align*}
\int \abs{\pi_n(x)}^2\, d\mu &\leq 2c\, \biggl\{\f{1}{n^{1+\alpha}} + \sum_{j=1}^{n-1}
\biggl( \f{j+1}{n}\biggr)^\alpha \biggl[ \f{1}{2n^2} + \f{1}{2j^2}\biggr]
\f{1}{n}\biggr\} \\
&\leq 2c\, \biggl\{\f{1}{2n^2} + \f{1}{n^{1+\alpha}}
\biggl( 1+\sum_{j=1}^\infty \f{(j+1)^\alpha}{j^2}\biggr)\biggr\} \\
&\leq C_1^{-1} n^{-(1+\alpha)}
\end{align*}
since $n^2\geq n^{1+\alpha}$ and $\sum_{j=1}^\infty \f{(j+1)^\alpha}{j^2} <\infty$
since $\alpha <1$.

By $\deg(\pi_n)=2n-2$, and Theorem~\ref{T.A.1},
\[
K_{2n-2} (0,0)^{-1} \leq C_1^{-1} n^{-(1+\alpha)}
\]
which is \eqref{7.26}.
\end{proof}

\section{Upper Bounds via OP Methods} \lb{s8}

Our main purpose in this section is to note that the upper bounds produced by
the method of Erd\"os-Turan \cite{ET40} provide a universal bound. So long as
$d\mu_\s =0$ near $x_0\in \supp(d\mu)$ and $w(x)$ is continuous and nonvanishing
at $x_0$, the bound is independent of the value of $w$ at $x_0$! Upper
bounds on spacing imply lower bounds on the density of zeros. Deift-Simon \cite{S169}
obtained universal lower bounds on the density of zeros, so the bounds we find are a
kind of microscopic analog of theirs. One key to the Erd\"os-Turan method is

\begin{lemma}\lb{L8.1} Let $x_1 < x_2 < \cdots < x_n$ in $\bbR$ and let $1\leq j\leq
n-1$. Then there exists a polynomial $\pi$ of degree at most $n-1$ so
\begin{align}
\pi(x_\ell) &= 0 \qquad  1\leq \ell\leq n; \, \ell\neq j, j+1 \lb{8.1} \\
\pi(x_j) &= \pi(x_{j+1}) =1 \lb{8.2} \\
\pi(y) & \geq 1 \quad\text{in}\quad [x_j, x_{j+1}] \lb{8.3}
\end{align}
\end{lemma}

\begin{proof} Let
\begin{equation} \lb{8.4}
\pi_0(x) =\prod_{\ell\neq j,j+1}\, (x-x_j)
\end{equation}
If
\begin{equation} \lb{8.5}
\pi_0 (x_j) =\pi_0 (x_{j+1})
\end{equation}
take
\begin{equation} \lb{8.5a}
\pi(x) = \f{\pi_0(x)}{\pi_0(x_j)}
\end{equation}
so \eqref{8.1}/\eqref{8.2} hold (we will look at \eqref{8.3} shortly). If \eqref{8.5} fails,
for $y\in\bbR\backslash [x_j, x_{j+1}]$, let
\begin{equation} \lb{8.6}
\pi_y(x) =(x-y) \pi_0(x)
\end{equation}
As $y$ runs through $(-\infty, x_j]$, $\f{y-x_j}{y-x_{j+1}}$ runs
from $1$ down to $0$, and as $y$ runs from $\infty$ to $x_{j+1}$,
the ratio runs from $1$ to $\infty$. Since \eqref{8.5} fails and
$\pi_0(x_j)$ and $\pi_0 (x_{j+1})$ have the same sign, there is a
unique $y$ with
\begin{equation} \lb{8.7}
\pi_y (x_j) =\pi_y(x_{j+1})
\end{equation}
so take
\begin{equation} \lb{8.8}
\pi(x) = \f{\pi_y(x)}{\pi_y(x_j)}
\end{equation}

In any event, $\pi$ obeys \eqref{8.1} and \eqref{8.2}. By Snell's theorem, $\pi'$ has a zero
between any two zeros of $\pi$, and so by counting degrees, exactly $1$. It follows that $\pi'$
has a local maximum in $[x_j, x_{j+1}]$ and no local minimum, so \eqref{8.3} holds.
\end{proof}

Let $\{x_j\}_{j=1}^n$ be the zeros of the OP, $P_n$, associated to a measure $d\mu$. Recall
(see Theorem~\ref{T.A.3}) that there are positive weights $\{\lambda_j\}_{j=1}^n$ so
\begin{equation} \lb{8.9}
\int \ti\pi (x)\, d\mu(x) = \sum_{\ell =1}^n \lambda_\ell \ti\pi (x_\ell)
\end{equation}
for any polynomial $\ti\pi$ with $\deg \ti\pi\leq 2n-1$.

\begin{theorem}\lb{T8.2} For any $j=1,2,\dots, n-1$,
\begin{equation} \lb{8.10}
\mu ([x_j, x_{j+1}]) \leq \lambda_j + \lambda_{j+1}
\end{equation}
\end{theorem}

\begin{proof} Let $\pi$ be the polynomial of degree $n-1$ or less given by Lemma~\ref{L8.1}.
Let $\ti\pi =\pi^2$, so $\deg\ti\pi \leq 2n-2\leq 2n-1$. Since $\ti\pi (x_\ell)=0$, $\ell\neq j,
j+1$, and $\ti\pi(x_j) =\ti\pi(x_{j+1})=1$,
\begin{equation} \lb{8.11}
\text{RHS of \eqref{8.9}}=\lambda _j + \lambda_{j+1}
\end{equation}
Since $\ti\pi\geq 0$ and $\ti\pi \geq 1$ on $[x_j,x_{j+1}]$,
\begin{equation} \lb{8.12}
\text{LHS of \eqref{8.9}} \geq \mu ([x_j, x_{j+1}])
\end{equation}
so \eqref{8.9} implies \eqref{8.10}.
\end{proof}

To exploit \eqref{8.10}, we need upper bounds on $\lambda_j$.

Suppose
\begin{equation} \lb{8.13}
E_\pm = {}_\text{\rm{inf}}^\text{\rm{sup}}\, \supp(d\mu )
\end{equation}
and for $E\in\ [E_-, E_+]$,
\begin{equation} \lb{8.14}
d(E)=\max (E_+-E, E-E_-)\leq E_+ -E_-
\end{equation}

\begin{theorem} \lb{T8.3} Suppose $I$ is a closed interval on which $d\mu$
is purely a.c.\ and $\max_{x\in I} w(x)=w_+<\infty$. Then for each $\delta >0$
and all weights $\lambda_j$ associated with $x_j\in I$ and $\dist
(x_j,\bbR\backslash I)\geq\delta$, we have with $m=[\f{n}{2}]$ and $n$
the number of zeros,
\begin{equation} \lb{8.15}
\lambda_j \leq \f{w_+ \pi \max_I d(E)}{m} + O\biggl( \f{1}{n^2}\biggr)
\end{equation}
where the $O(\f{1}{n^2})$ is uniform in all $\lambda$'s with the given $\delta$
{\rm{(}}and depends on $\max_I d(E)${\rm{)}}.
\end{theorem}

\begin{proof} Let $a=\max_I d(E)$ and $x_0$ be the $x_j$ for $\lambda_j$, and let
$\ti\pi(x)= \pi_m (x;x_0,a)$ given by Theorem~\ref{T.A.6}. Since $\ti\pi(x)=1$ and
$\deg (\ti\pi^2)\leq 2n-1$, we have
\begin{align*}
\lambda_j &\leq \sum \lambda_\ell \ti\pi(x_\ell)^2 \\
&=\int d\mu(x) \ti\pi (x)^2 = K_1 + K_2
\end{align*}
where $K_1$ is the integral over $(x_0-\delta, x_0+\delta)$, and $K_2$ the integral over all
other $x$. By \eqref{A.1.33},
\[
K_2 = O\biggl( \f{1}{n^2}\biggr)
\]
with estimates only dependent on $a$ and $\delta$. For $\ti\pi(x)^2=O(\f{1}{n^2})$ on the region of
integration and $\mu(\bbR)=1$.

If $(x_0-\delta, x_0+\delta)\subset I$, $d\mu\leq w_+\, dx$, so \eqref{8.15} follows from
\eqref{A.1.34}.
\end{proof}

\begin{theorem}\lb{T8.4} Suppose $I$ is a closed interval on which $d\mu$ is purely a.c.\ and
\begin{equation} \lb{8.16}
0 < w_- \equiv\min_{x\in I}\, w(x) \leq \max_{x\in I}\, w(x)\equiv w_+ < \infty
\end{equation}
Then for any $E\in I^\intt$,
\begin{equation} \lb{8.17}
\limsup_{n\to\infty}\, n[z_n^{(1)}(E) - z_n^{(-1)} (E)]\leq 4\pi \, d(E) \lim_{\delta\downarrow 0}\,
\f{\max\{w(x)\mid\abs{x-E}<\delta\}}{\min\{w(x)\mid\abs{x-E}<\delta\}}
\end{equation}
In particular, if $E$ is a point of continuity of $w$,
\[
\text{LHS of \eqref{8.17}}\leq 4\pi\, d(E)
\]
independently of the value of $w(E)$.
\end{theorem}

\begin{proof} Clearly,
\begin{equation} \lb{8.18}
\begin{split}
\abs{z_n^{(1)}(E) & -z_n^{(-1)}(E)} \\
& \leq [\min\{w(x)\mid
z_n^{(-1)}(E)\leq x \leq z_n^{(1)}(E)\}]^{-1} \mu(z_n^{(-1)}(E),
z_n^{(1)}(E))
\end{split}
\end{equation}
From this, \eqref{8.10}, \eqref{8.15}, and $\lim n/m=2$, we get \eqref{8.17} by using the fact
that since $E$ is a limit point of an infinity of zeros, we have $\lim_{n\to\infty}
\abs{z_n^{\pm 1} (E)-E} =0$.
\end{proof}

This is the promised universal lower bound on the density of zeros. The method is flexible enough
to say something if $w(x)$ has a zero of a fixed order.

\begin{theorem} \lb{T8.5} Suppose $d\mu$ is purely absolutely continuous in a neighborhood
of $E_0$, and for some $q >0$,
\begin{equation} \lb{8.19}
0<\gamma_- =\liminf_{x\to E_0}\, \f{w(x)}{\abs{x-E_0}^q} \leq \limsup_{x\to E_0}\,
\f{w(x)}{\abs{x-E_0}^q} = \gamma_+ < \infty
\end{equation}
Then
\begin{equation} \lb{8.20}
\limsup\, n\abs{z_n^{(1)}(E_0) -z_n^{(-1)}(E_0)} <\infty
\end{equation}
\end{theorem}

\begin{proof} By \eqref{8.19} for any $\delta$, there is $N$ so for $n\geq N$,
\begin{equation} \lb{8.21}
\mu([z_n^{(-1)}(E_0), z_n^{(1)}(E_0)]) \geq (\gamma_- -\delta) 2^{-q}(q+1)^{-1} (z_n^{(1)}(E_0)
-z_n^{(-1)}(E_0))^{q+1}
\end{equation}
By the proof of \eqref{8.15}, the $\lambda$'s associated to $z_n^\pm (E)$ obey
\begin{equation} \lb{8.22x}
\lambda_j \leq \f{C_1 (\gamma_+ + \delta)[\max (C_2 n^{-1}, \abs{z_n^\pm (E_0)-E_0})]^q}{n}
\end{equation}
for constants $C_1, C_2$. \eqref{8.10},\eqref{8.21}, and \eqref{8.22x} imply \eqref{8.20}.
\end{proof}

Given our bounds in Section~\ref{s7}, we can also say something when the singularity of
the weight is not as regular as some power. The key is an abstraction of an argument of
Nevai \cite{Nev79} (see also Golinskii \cite{Gol02}).

\begin{theorem}\lb{T8.6} Let $a=\max (\supp(d\mu))-\min (\supp(d\mu))$. Fix integers
$p,q$ so that
\begin{equation} \lb{8.22}
(2p-2)^{2q} \leq 2n-1
\end{equation}
Then for any successive zeros $E,E'$ of $p_n$, we have
\begin{equation} \lb{8.23}
\abs{E-E'} \leq \f{a}{p}\, [K_p (\tfrac12\, (E+E'), \tfrac12\,(E+E'))]^{1/2q}
\end{equation}
\end{theorem}

\begin{proof} Let $\ti\pi$ be defined in terms of the $\pi$ of Theorem~\ref{T.A.6} by
Theorem~\ref{T.A.3},
\begin{equation} \lb{8.24}
\ti\pi (x)=[\pi_p (x;\tfrac12\, (E+E'),a)]^q
\end{equation}
By \eqref{8.22}, $\deg [\ti\pi]^2\leq 2n-1$ so, by \eqref{8.9},
\begin{align}
\int \abs{\ti\pi(x)}^2\, d\mu(x) &=\sum_{j=1}^n \lambda_j \abs{\ti\pi (E_j)}^2 \lb{8.25}  \\
&\leq \biggl( \f{1}{2p} + \f{a}{2p \abs{E-E'}}\biggr)^{2q} \lb{8.26}
\end{align}
since $\sum\lambda_j=1$, $\lambda_j\geq 0$, and $\min \abs{E_j -\f12 (E+E')}=\f12
\abs{(E-E')}$.

Since $\ti\pi (\f12 (E+E'))=1$, by Theorem~\ref{T.A.1},
\begin{equation} \lb{8.27}
K_n (\tfrac12\, (E+E'), \tfrac12\, (E+E'))^{-1} \leq \int \abs{\ti\pi(x)}^2\, d\mu(x)
\end{equation}

Since $\abs{E-E'}\leq a$,
\begin{equation} \lb{8.28}
\f{1}{2p} + \f{a}{2p\abs{E-E'}} \leq \f{a}{p\abs{E-E'}}
\end{equation}
\eqref{8.23} is immediate from \eqref{8.26}, \eqref{8.27}, and \eqref{8.28}.
\end{proof}

The following abstracts an argument of Golinskii, who needed to make global hypotheses since
he relied on estimates of Geronimus:

\begin{corollary}\lb{C8.7} Suppose that for some interval $I$, $A>0$ and $C$, we have
\[
\sup_{E\in I}\, \abs{K_n (E,E)}\leq C(n+1)^A
\]
Then for any $\delta >0$,
\begin{equation} \lb{8.29}
\lim_{n\to\infty}\, \sup_{\substack{ E,E'\text{ successive zeros} \\
\dist(E,\bbR\backslash I) >\delta }} \abs{E-E'} \biggl[ \f{n}{\log n}\biggr] <\infty
\end{equation}
\end{corollary}

\begin{proof} Pick $q =[\log n]$ and $p$ as large as possible so that \eqref{8.22} holds.
Since
\[
(n+1)^A = \exp (A\log (n+1))
\]
$[K_p]^{1/2q}$ is bounded and \eqref{8.23} implies \eqref{8.29}.
\end{proof}

Combining this corollary and Theorem~\ref{T7.1}, we obtain a local version of
Golinskii's \cite{Gol02} result:

\begin{corollary}\lb{C8.8} If \eqref{7.1} and \eqref{7.2} hold, then we have \eqref{8.29}
for $I=(x_0-a, x_0+a)$.
\end{corollary}

We also have the following (a local version of results of Nevai
\cite{Nev79} and Golinskii \cite{Gol02}):

\begin{theorem}\lb{T8.9} Suppose for some interval $I$ we have that
\begin{equation} \lb{8.30}
d\mu =w\, dx +d\mu_\s
\end{equation}
where
\begin{equation} \lb{8.31}
\int_I \log w\, dx >-\infty
\end{equation}
Then for any $\delta >0$,
\begin{equation} \lb{8.31b}
\lim_{n\to\infty}\, \sup_{\substack{ E,E'\text{ successive zeros} \\
\dist (E,\bbR\backslash I)>\delta }} \abs{E-E'} n^{1/2} <\infty
\end{equation}
\end{theorem}

We need the following lemma:

\begin{lemma}\lb{L8.10} If $J$ has Jacobi parameters obeying
\begin{equation} \lb{8.31a}
\sum_{n=1}^\infty\, \abs{a_n-1}^2 + \abs{b_n}^2 <\infty
\end{equation}
then for any $\delta >0$, there is $C_\delta$ so
\begin{equation} \lb{8.32}
\sup_{E\in [-2+\delta, 2-\delta]}\, \|T_n(E)\|\leq \exp \bigl( C_\delta
\sqrt{n+1}\,\bigr)
\end{equation}
\end{lemma}

\begin{proof} Define $u_n^\pm (E)=e^{in\theta}$ where $2\cos\theta = E$ and $0<\theta
<\pi$. By standard variation of parameters about $u_n^\pm$, one proves
\[
\text{LHS of \eqref{8.32}} \leq \prod_{j=1}^n \{1+C(\abs{b_j}+\abs{a_j-1})\}
\]
Since $1+x \leq \exp(x)$ and
\[
\sum_{j=1}^n (\abs{b_j} + \abs{a_j-1}) \leq \biggl[\, \sum_{j=1}^n \, (\abs{b_j}^2 +
\abs{a_j-1}^2)\biggr]^{1/2} [2n]^{1/2}
\]
\eqref{8.32} is immediate.
\end{proof}

\begin{proof}[Proof of Theorem~\ref{T8.9}] By scaling, we suppose $I=[-2,2]$. By
Corollary~\ref{C.A.2},
\begin{equation} \lb{8.33}
K_n(x,x;d\mu) \leq K_n (x,x;d\nu)
\end{equation}
where
\begin{equation} \lb{8.34}
d\nu = \chi_{[-2,2]} w\, dx
\end{equation}

Let $d\ti\nu$ be the normalized $d\nu$. By the theorem of Killip-Simon \cite{KS},
the Jacobi parameters obey \eqref{8.31a}, so by \eqref{8.33},
\[
\sup_{x\in [-2+\delta, 2-\delta]}\, K_n (x,x;d\mu)\leq \exp\bigl( C_\delta
\sqrt{n+1}\,\bigr)
\]

In Theorem~\ref{T8.6}, take $q=\bigl[\sqrt{n}\bigr]$ and $p$ as large as can be so
\eqref{8.22} holds. Then $p\sim c\sqrt{n}$ and \eqref{8.23} implies \eqref{8.31b}.
\end{proof}

\section{Lower Bounds via OP Methods} \lb{s9}

In this section, we will get lower bounds in terms of the CD kernel alone. The basic
method is due to Golinskii \cite{Gol02}, but when he applied the method, he made global
assumptions on the measure, and we want to note that local assumptions suffice. Other OP
lower bound methods are due to Erd\"os-Turan \cite{ET40} and Nevai \cite{Nev79}.

\begin{theorem}\lb{T9.1} If $E,E'$ are distinct zeros of $P_n(x)$, if $\bar E=\f12 (E+E')$
and $\delta >\f12 \abs{E-E'}$, then
\begin{equation} \lb{9.1}
\abs{E-E'} \geq \f{[\delta^2 - (\f12 \abs{E-E'})^2]}{3n}\, \biggl[ \f{K_n (E,E)}
{\sup_{\abs{y-\bar E}\leq \delta} K_n (y,y)}\biggr]^{1/2}
\end{equation}
\end{theorem}

\begin{remarks} 1. In most applications, $\delta$ is fixed and $\abs{E-E'} \to 0$, so
$\delta^2 -(\f12 \abs{E-E'})^2\sim \delta^2 >0$. In typical cases, the $\inf$ and
$\sup$ of $K_n (y,y)$ for $\abs{y-\bar E}<\delta$ are comparable and \eqref{9.1} gives
an $1/n$ lower bound.

2. This theorem also yields a result with the same asymptotics as \eqref{2.22} for
$K_n(E,E)\geq 1$, while the $\sup$ is bounded exponentially in $n$.

3. Interestingly enough, the proof here depends on \eqref{2.19} written as \eqref{9.2}.

4. It is interesting to compare Theorems~\ref{T9.1} and \ref{T2.2}. \eqref{2.17} only
depends on information at $E_0$ while \eqref{9.1} has a $\sup K_n (y,y)$ over a neighborhood,
but \eqref{2.17} requires information on both solutions of \eqref{2.1} while \eqref{9.1} only
on $p_n$.
\end{remarks}

\begin{proof} Since $p_n(E)=p_n (E')=0$ and $E\neq E'$, we have
\begin{equation} \lb{9.2}
K_n (E,E')=0
\end{equation}
by \eqref{A.1.3} (the Christoffel-Darboux formula). Thus, (supposing $E<E'$ for notational
convenience),
\begin{align}
K_n (E,E) &= K_n(E,E) - K_n (E,E')  \lb{9.3} \\
&\leq \abs{E-E'} \sup_{\abs{y-\bar E}\leq \f12\abs{E-E'}} \biggl| \f{\partial}{\partial y}\,
K_n (E,y)\biggr| \lb{9.4} \\
&\leq \abs{E-E'} \{\delta^2 -[\tfrac12\, (E-E')]^2\}^{-1} \notag \\
& \qquad\qquad \sup_{\abs{y-\bar E}\leq \delta}\, [\delta^2 - (y-\bar E)^2 ]^{1/2}
\biggl| \f{\partial}{\partial y}\, K_n(E,y)\biggr| \lb{9.5}
\end{align}

$K_n (E,y)$ is a polynomial in $y$ of degree $n$, so by
\eqref{A.1.28} and \eqref{9.5},
\begin{equation} \lb{9.6}
K_n (E,E)\leq \abs{E-E'} \{\delta^2 - [\tfrac12\, (E-E')]^2\}^{-1}
(3n) \sup_{\abs{y-\bar E}<\delta}\, \abs{K_n (E,y)}
\end{equation}

By the Schwartz inequality,
\begin{equation} \lb{9.7}
\sup_{\abs{y-\bar E}<\delta }\, \abs{K_n (E,y)} \leq K_n (E,E)^{1/2}
\sup_{\abs{y-\bar E}<\delta}\, \abs{K_n (y,y)}^{1/2}
\end{equation}
\eqref{9.6} plus \eqref{9.7} imply \eqref{9.1}.
\end{proof}

\begin{corollary} \lb{C9.2} Let $I$ be an interval on which
\begin{equation} \lb{9.8}
t=\sup_{n,y\in I}\, \|T_n(y)\| <\infty
\end{equation}
Then, for any $E\in I^\intt$,
\begin{equation} \lb{9.9}
\liminf_{n\to\infty} \, n\abs{z_n^{(+1)}(E) -z_n^{(-1)}(E)} \geq
\f{t^{-2}}{3}\, \dist (E,\bbR\backslash I)^2
\end{equation}
\end{corollary}

\begin{remark} This should be compared with what follows from \eqref{2.15} and \eqref{2.17}
which implies
\[
\text{LHS of \eqref{9.9}} \geq t^{-2}
\]
\end{remark}

\begin{proof} \eqref{9.9} follows from \eqref{9.1} if one notes that for $y\in I$,
\[
(n+1) t^{-2} \leq K_n (y,y)\leq (n+1)t^2
\qedhere
\]
\end{proof}

We can also use Theorem~\ref{T9.1} to get a lower bound in terms of local bounds on
the weights.

\begin{theorem}\lb{T9.3} Suppose $d\mu =w\, dx + d\mu_\s$ where $d\mu_\s (x_0-\delta,
x_0+\delta)=0$ and
\begin{equation} \lb{9.10}
0 < \inf_{\abs{y-x_0}\leq \delta}\, w(x) \leq \sup_{\abs{y-x_0}\leq\delta}\, w(x) <\infty
\end{equation}
Then for any $\veps <\delta$,
\begin{equation} \lb{9.11}
\inf_{\abs{y-x_0}<\veps}\, \liminf_{n\to\infty} \, n\abs{z_n^{(+1)}(y) - z_n^{(-1)}(y)} >0
\end{equation}
\end{theorem}

\begin{proof} By \eqref{9.1}, it suffices to prove
\begin{equation} \lb{9.12}
\sup_{\abs{y-x_0}<\veps}\, [n^{-1} K_n(y,y)] <\infty
\end{equation}
and
\begin{equation} \lb{9.13}
\inf_{\abs{y-x_0}<\veps}\, [n^{-1} K_n (y,y)] >0
\end{equation}

By \eqref{9.10}, for $\veps$ fixed, uniformly in $y$ with $\abs{y-x_0}<\veps$, we can
find a fixed scaling and some translate of $c(4-x^2)^{1/2} \chi_{[-2,2]}\, dx$ lying
below $d\mu (x-y)$, so using \eqref{A.1.13} and the explicit $K_n$ for Chebyshev
polynomials of the second kind (i.e., the free $K_n$), we get \eqref{9.12}.

On the other hand, by \eqref{A.1.4}/\eqref{A.1.5},
\begin{equation} \lb{9.14}
K_n (x_0, x_0) \geq \biggl( \int \abs{\pi_n(x)}^2\, d\mu(x)\biggr)^{-1}
\end{equation}
for any $\pi_n$ of degree $n$ with $\pi_n (x_0)=1$. Using a suitable $\pi_{n/2}$
of the form given by Theorem~\ref{T.A.6} and estimates we used earlier in this
paper, we get an  $O(n)$ lower bound on $K_n$, that is, \eqref{9.13} holds.
\end{proof}

\section{Zeros of POPUC} \lb{s10}

While we have discussed OPRL up to now, virtually all the ideas extend to POPUC. POPUC
are defined by taking the first $n-1$ recursion parameters (Verblunsky coefficients),
$\alpha_0, \dots, \alpha_{n-2}$ in $\bbD$, and picking $\beta\in\partial\bbD$ and letting
\begin{equation} \lb{10.1}
\ti\varphi_n(z)=z\varphi_{n-1}(z) -\bar\beta \varphi_{n-1}^*(z)
\end{equation}
Since $\varphi_{n-1}^*$ is nonvanishing on $\bbD$ and
$\abs{\varphi_{n-1}}= \abs{\varphi_{n-1}^*}$ on $\partial\bbD$, we
have $\abs{\varphi_{n-1}/\varphi_{n-1}^*} <1$ on $\bbD$ by the
maximum principle. Thus, $\ti\varphi_n(z)$ is nonvanishing on $\bbD$
and, by symmetry ($\ol{\ti\varphi_n(1/\bar z)}=z^{-n}
[\varphi_{n-1}^*-\beta z \varphi_{n-1}]$), nonvanishing on
$\bbC\backslash\ol{\bbD}$. Thus, the zeros of $\ti\varphi_n$ lie on
$\partial\bbD$; indeed, they are the zeros of a finite unitary
matrix (see Theorem~8.2.7 of \cite{OPUC1}). Zeros of POPUC are
discussed extensively in Golinskii \cite{Gol02},
Cantero-Moral-Vel\'azquez \cite{CMV02,CMV06}, Simon \cite{Roppopuc},
and Wong \cite{Wong}.

As explained in Section~10.8 of \cite{OPUC2}, there is an OPUC analog of
\eqref{2.7}--\eqref{2.8} (namely, (10.8.3)--(10.8.5) of \cite{OPUC2}) which immediately
leads to an analog of Theorem~\ref{T2.1}. While (10.8.3)--(10.8.5) are stated for
the solutions $\psi_\bddot + F\varphi_\bddot$, they also hold for $\psi_\bddot$ and
$\varphi_\bddot$. Key to this analog is the orthogonality of $\varphi_\bddot (z_1)$
and $\varphi_\bddot (z_2)$ for two zeros of $\ti\varphi_n(z)$. This follows from the
CD formula for OPUC (see Theorem~2.2.7 of \cite{OPUC1}) for $\ti\varphi_n(z_1) =
\ti\varphi_n(z_2)=0$ implies $z_j\varphi_n(z_j)=\bar\beta\varphi_n^*(z_j)$, so
\[
\ol{\varphi_n^*(z_1)}\, \varphi_n^*(z_2) -\bar z_1 z_2\, \ol{\varphi_n(z_1)}\,
\varphi_n(z_2) =0
\]
and thus, by (2.2.42) of \cite{OPUC1},
\begin{equation} \lb{10.2}
\sum_{j=0}^n\, \ol{\varphi_j (z_1)}\, \varphi_j (z_2) =0
\end{equation}

Combined with the techniques of Section~12.1 of \cite{OPUC2} and our proof of
Theorem~\ref{T3.4}, we get

\begin{theorem}\lb{T10.1} Let $\{\alpha_n\}_{n=0}^\infty$ be a set of Verblunsky
coefficients that obeys
\begin{equation} \lb{10.3}
\sum_{n=0}^\infty\, \abs{\alpha_{n+1}-\alpha_n} <\infty
\end{equation}
and
\begin{equation} \lb{10.4}
\abs{\alpha_n}\to 0
\end{equation}
Then the zeros of the POPUC, $\ti\varphi_n(z)$, for any choice of $\beta$ have uniform
clock behavior on any compact subset of $\partial\bbD\backslash\{1\}$.
\end{theorem}

\begin{remarks} 1. An interesting example is $\alpha_n =(n+2)^{-\beta}$ for any $\beta >0$.
This is related to a conjecture of \cite{Saff1}, albeit the conjecture there is for OPUC,
not POPUC.

2. The density of zeros in this case is $d\theta/2\pi$ on $\partial\bbD$; see Theorem~8.2.7
and Example~8.2.8 of \cite{OPUC1}.

3. If \eqref{10.3}/\eqref{10.4} are replaced by
\begin{equation} \lb{10.5}
\sum_{n=0}^\infty \, \abs{\alpha_n}<\infty
\end{equation}
then $\partial\bbD\backslash\{1\}$ can be replaced by $\partial\bbD$. This is a result of
\cite{Saff1}. Because we have global control in this case, one does not need a priori $1/n$
bounds on zero spacing.
\end{remarks}

There are also analogs of the bounds of Sections~\ref{s5}, \ref{s6}, \ref{s8}, and \ref{s9}:

\smallskip
1.  One has that
\begin{equation} \lb{10.6}
\f{d^2}{d\theta^2}\, \log \abs{e^{i\varphi}-e^{i\theta}} =
-\f{1}{\abs{e^{i\varphi}- e^{i\theta}}^2}
\end{equation}
so there is a bound like \eqref{5.4} for POPUC (all of whose zeros lie on $\partial\bbD)$, and
thus, there is an analog of Theorem~\ref{T5.1}.

\smallskip
2. If one defines $\eta_n(z)$ by
\begin{equation} \lb{10.7}
e^{i\eta_n(\theta)} = \f{e^{i\theta}\varphi_{n-1} (e^{i\theta})}
{\varphi_{n-1}^* (e^{i\theta})}
\end{equation}
then
\begin{equation} \lb{10.8}
\f{d\eta_n}{d\theta} = \f{[\sum_{j=0}^{n-1} \abs{\varphi_j (e^{i\theta})}^2]}
{\abs{\varphi_{n-1} (e^{i\theta})}^2}
\end{equation}
This follows from (2.2.71) of \cite{OPUC1} which implies
\begin{equation} \lb{10.9}
\left. \f{\partial}{\partial r}\, \log \abs{\varphi_{n+1}(re^{i\theta})}^2\right|_{r=1}
= (n+1) + \abs{\varphi_{n+1} (e^{i\theta})}^{-2} \sum_{j=0}^n\, \abs{\varphi_j (e^{i\theta})}^2
\end{equation}

By the Cauchy-Riemann equations
\begin{equation} \lb{10.10}
\text{LHS of \eqref{10.9}} = 2\, \f{\partial}{\partial\theta}\,
\arg [\varphi_{n+1} (e^{i\theta})]
\end{equation}
Since \eqref{10.7} implies
\begin{equation} \lb{10.11}
\eta_n = \theta - (n-1)\theta + 2 \arg [\varphi_{n-1}(e^{i\theta})]
\end{equation}
we obtain \eqref{10.8}.

\eqref{10.8} implies $d\eta_n/d\theta >0$ and, given that zeros of $\ti\varphi_n$ occur
when $\eta_n =\arg \bar\beta$ (mod $2\pi$), bounds like those of Theorem~\ref{T6.3} on
zero spacing for POPUCs.

\smallskip
3. Since all the techniques of the Appendix extend to OPUC, the estimates of
Sections~\ref{s8} and \ref{s9} extend to POPUC; indeed, somewhat weaker variants
occur already in Golinskii \cite{Gol02}.

\section{Examples, Counterexamples, Conjectures, and Questions} \lb{s11}

\noindent (a) {\it Clock behavior based only on local behavior}. \ Suppose
\eqref{7.1} holds on for some $(c,d)\subset\bbR$, we have $d\mu_\s ([c,d])=0$,
$w>0$ on $(c,d)$ and $w$ is $C^\infty$ there. We have proven $O(1/n)$ upper and
lower bounds in this case. This leads to the natural question:

\begin{oq}\lb{Q11.1} Under the above hypothesis, does one have clock behavior on $(c,d)$?
\end{oq}

This is a very subtle question because clock behavior involves the density of
states, and it is not even clear that exists on $(c,d)$ only under the above
hypothesis. What is clear is that if the density of states exists, it is a
global quantity and not just dependent on $w$ on $(c,d)$. We want to
demonstrate this by example. We will need the following:

\begin{proposition}\lb{P11.2} Let $d\mu_0$ be given by
\begin{equation} \lb{11.1}
d\mu_0 = (2\pi)^{-1} \sqrt{4-x^2}\, \chi_{[-2,2]}\, dx
\end{equation}
Suppose $f$ is a $C^2$ function on $[-2,2]$ with $f'(2) = f'(-2) =0$, $f''(-2)
= f''(2)$, and
\begin{equation} \lb{11.2}
f\geq \alpha > 0 \qquad \int f\,d\mu_0 =1
\end{equation}
where $\alpha >0$ is a positive real. Let
\begin{equation} \lb{11.3}
d\mu(x) = f(x)\, d\mu_0 (x)
\end{equation}
Then the density of states exists for $d\mu$ and is given by \eqref{1.7} and
there is clock behavior uniformly on each interval $[-2+\veps, 2-\veps]$.
\end{proposition}

\begin{proof} By Theorem~13.2.1 of \cite{OPUC2}, there is a map $\Sz_2$
of real measures, $d\rho$, on $\partial\bbD$ (i.e., those measures with real
Verblunsky coefficients) to those measures $d\mu$ on $[-2,2]$ which are of the
form $f\, d\mu_0$ where
\begin{equation} \lb{11.4}
\int_{-2}^2 f(x) (4-x^2)^{-1/2}\, dx <\infty
\end{equation}
and under this map, $d\rho(\theta) = w(\theta)\, d\theta$ where
\begin{equation} \lb{11.5}
w(\theta) = c\, f(2\cos\theta)
\end{equation}

Thus, $w$ is $C^2$ (including at $\theta=0,\pi)$ and so
$\sum_{n=-\infty}^\infty n^2 \abs{\widehat w_n}^2 <\infty$, which implies
$\sum_{n=-\infty}^\infty \abs{\widehat w_n}<\infty$. Thus, by Baxter's theorem
(Theorem~5.2.1 of \cite{OPUC1}), the Verblunsky coefficients are in $\ell^1$.
By (13.2.20)/(13.2.21) of \cite{OPUC2}, \eqref{1.14} holds, which implies the
claimed result.
\end{proof}

\begin{example}\lb{E11.3} Let $d\ti\mu_0$ be $d\mu_0$ scaled to $[-1,1]$. We
pick $f_1$ obeying the hypothesis of Proposition~\ref{P11.2} and $f_2$ scaled
to $[-1,1]$, and so that $d\mu_1 =f_1\, d\mu_0$ and $d\mu_2 = f_2\, d\ti\mu_0$
obey $d\mu_1 \equiv d\mu_2$ on $[-\f12,\f12]$. Both have clock behavior on
$[-\f12, \f12]$ but with different density of states, namely \eqref{1.7} and
\eqref{1.7} scaled. \qed
\end{example}

\medskip
\noindent (b) {\it Pointwise upper bounds}. \ We obtained lower bounds on
$z_n^{(1)} (E_0)-z_n^{(-1)}(E_0)$ if $T_n (E_0)$ is bounded, but our upper
bounds required control of $T_n(E)$ for $E$ in a neighborhood of $E_0$.

\begin{oq}\lb{Q11.4} Are there upper bounds on spacing if we only know that
$T_n (E_0)$ is bounded?
\end{oq}

\medskip
\noindent (c) {\it Improved spacing estimates}.

\begin{oq}\lb{Q11.5} Can $\sqrt{n}$ in \eqref{2.25} be improved?
\end{oq}

\medskip
\noindent (d) {\it More on spacing and $\gamma(E)$}. \ We saw that $\gamma
(E_0) =0$ plus some regularity of $\nu$ at $E_0$ implies an $O(1/n)$ lower
bound. Does it imply clock spacing? In particular,

\begin{oq}\lb{Q11.6} Is there local clock behavior for a.e.\ $E_0$ with $\gamma
(E_0)=0$ in the case of almost periodic Jacobi parameters?
\end{oq}

\begin{example}\lb{E11.7} \cite{KStoi} has proven, for $\alpha <\f12$, the
OPUC analog of Poisson behavior for $a_n\equiv 1, b_n$ independent random
variable of the form $b_n =Cn^{-\alpha} w_n$ where $w_n$ is uniformly
distributed in $[-1,1]$. We assume their result is true in the Jacobi case. Of
course, $\gamma (E)=0$ in this case. We do not have a contradiction with
Corollary~\ref{C4.2} since this model is not ergodic. The example does show
though that ergodicity is a necessary hypothesis. \qed
\end{example}


\begin{example}\lb{E11.7a} Corollary \ref{C4.2} shows that ergodicity
along with Poisson behavior imply positive Lyapunov exponent. This raises the
natural question: Does ergodicity along with a positive Lyapunov exponent
imply Poisson behavior? The answer is negative, as can be shown by the
following example: Consider the Jacobi matrix with $a_n=1$,
$b_n = \lambda\cos(2\pi\alpha n + \theta)$, where $|\lambda| > 2$ and $\alpha$
is a (Liouville) irrational for which there is a sequence of rationals
$\{p_n/q_n\}_{n=1}^\infty$ such that $|\alpha - p_n/q_n|<n^{-q_n}$. This is
an ergodic Jacobi matrix and it is well known (see, e.g., \cite{S149}) that
its Lyapunov exponent $\gamma(E)$ is positive for any $E$. By using the
results of Avron-van Mouche-Simon \cite{S211} and considering scales of the form
$mq_n$, where $m > 2$ is an integer, one can show that, for each $\theta$,
there would be clusters of $m-1$ zeros each of which is contained in
an interval whose length is of order $(2/|\lambda|)^{q_n/2}$. As $\theta$
is varied, these clusters will move over regions whose size is roughly
of order $1/q_n$. This behavior contradicts Poisson behavior. More precisely,
it is possible to show that Poisson behavior does not occur for Lebesgue
a.e.\ $E$ in the spectrum.
\end{example}

\medskip
\noindent (e) {\it Zero spacing and the Szeg\H{o} condition}.

\begin{oq}\lb{Q11.8} Does one have $O(n^{-1})$ bounds (upper and lower) when a
Szeg\H{o} or quasi-Szeg\H{o} condition holds?
\end{oq}

\medskip
\noindent (f) {\it Spacing at zeros of $w(x)$}.

\begin{oq}\lb{Q11.9} What can one say at zero spacing at points $x_0$ where
$w(x)$ has a ``regular" zero, that is, $w(x)\sim \abs{x-x_0}^\alpha$ for some
$\alpha >0$?
\end{oq}

\medskip
\noindent (g) {\it Edge zeros when $a_n = 1-n^{-\gamma}$}. \ The following
illuminates Theorem~\ref{T3.8}.

\begin{example}\lb{E11.10} Let $b_n\equiv 0$, $a_n = 1-n^{-\gamma}$ for $\gamma
>0$. Then Theorem~\ref{T3.4} applies and there is clock behavior away from $-2$
and $2$. If $\gamma >2$, Theorem~\ref{T3.8} applies and the largest $E_j$ has
$E_j =2-Cn^{-2} + o(n^{-2})$. We claim for $\gamma$ in general
\begin{equation} \lb{11.6}
2-C_2 n^{-\gamma} + O(n^{-1}) \leq E_j^{\max} \leq 2 - C_1 n^{-\gamma}
\end{equation}
capturing the leading behavior for $\gamma\in (0,1)$. The upper bound in
\eqref{11.6} comes from monotonicity of $E_j^{\max}$ in the $a$'s and the fact
that for $J_{n;F} \max_{1\leq j\leq n} \abs{a_j}\sim 1-n^{-\gamma}$. The lower
bound comes with a trial vector that lives in $[n/2,2]$ with maximum at $3n/4$
and constant slope in between. \qed
\end{example}

\medskip
\noindent (h) {\it OPUC}. \ This paper has a fairly complete analysis of OPRL
and POPUC. Many questions remain for general OPUC.

\appendix
\section{Tools of the OP Trade} \lb{App}
\renewcommand{\theequation}{A.\arabic{equation}}
\renewcommand{\thetheorem}{A.\arabic{theorem}}
\setcounter{theorem}{0}
\setcounter{equation}{0}

As explained in the introduction, this paper is intended for two audiences,
so we include here a summary of tools well known to the OP community but not
so well to the Schr\"odinger operator community. Because the tools, while powerful,
are simple, we can even give complete proofs. We will discuss the Christoffel
variational principle, Gauss-Jacobi quadrature, Bernstein's inequality, and
Dirichlet-Fej\'er trial polynomials.

\smallskip
\noindent {\bf (a) Christoffel variational principle}. \ We will define OPRL for
arbitrary positive measures (with finite moments) even if $\mu(\bbR)\neq 1$. The
monic polynomials, $P_n$, are independent of normalization, but the orthonormal
polynomials, $p_n$, are not. For example, $p_0(x,d\mu)=d\mu(\bbR)^{-1/2}$. The
{\it Christoffel-Darboux} (a.k.a.\ CD) {\it kernel\/} or {\it reproducing kernel\/}
is defined by
\begin{equation} \lb{A.1.1}
K_n (x,y)=\sum_{j=0}^n\, \ol{p_j(y)}\, p_j(x)
\end{equation}
We will use $K_n(x,y;d\mu)$ if the measure needs to be more explicit. The name
reproducing kernel comes from
\begin{equation} \lb{A.1.2}
(P_n f)(x) =\int K_n(x,y) f(y)\, d\mu(y)
\end{equation}
where $P_n$ is the projection in $L^2 (\bbR,d\mu)$ onto the space of polynomials
of degree $n$.

We need the following in Section~\ref{s6}:

\begin{theorem}[CD formula] \lb{T.A.1A} We have
\begin{equation} \lb{A.1.3}
K_n(x,y)=a_{n+1} \biggl[ \f{\ol{p_{n+1}(y)}\, p_n(x) - \ol{p_n(y)}\, p_{n+1}(x)}
{\bar y-x}\biggr]
\end{equation}
\end{theorem}

\begin{proof} This is a discrete version of integrating a Wronskian. Take the
equation \eqref{1.3} at $\bar y$, multiply by $p_n(x)$ and subtract \eqref{1.3}
at $x$, multiplied by $\ol{p_n(y)}$, and obtain
\[
Q_{n+1}(x,y) = \ol{p_n(y)}\, p_n(x) + Q_n (x,y)
\]
where $Q_{n+1}$ is the right side of \eqref{A.1.3}. Since $p_{-1}(x)=0$, $Q_0
\equiv 0$, so \eqref{A.1.3} follows by iteration.
\end{proof}

\begin{remark} For $x,y$ real, $p_n(x), p_n(y)$, are real so the bars are not needed.
Indeed, one can drop all the bars for complex $x,y$, but given the OPUC analogs, it is
natural to use the bars.
\end{remark}

\begin{theorem}[Christoffel variational principle] \lb{T.A.1} Let
\begin{equation} \lb{A.1.4}
\lambda_n (x_0;d\mu) =\inf\biggl( \int \abs{\pi_n(x)}^2\, d\mu(x) \biggm|
\deg\pi_n\leq n;\, \pi_n(x_0)=1\biggr)
\end{equation}
Then
\begin{equation} \lb{A.1.5}
\lambda_n (x_0;d\mu) =K_n (x_0,x_0;d\mu)^{-1}
\end{equation}
\end{theorem}

\begin{remarks} 1. $\lambda_n$ are called {\it Christoffel numbers}. More generally,
we have $p$-Christoffel numbers defined, for $0<p<\infty$, by
\begin{equation} \lb{A.1.6}
\lambda_n (x_0,p;d\mu) =\inf \biggl( \int \abs{\pi_n(x)}^p\, d\mu(x) \biggm|
\deg \pi_n \leq n; \, \pi_n(x_0)=1\biggr)
\end{equation}

2. Another way of writing \eqref{A.1.4} is that $\lambda_n$ is the optimal constant in
\begin{equation} \lb{A.1.7}
\abs{\pi_n(x_0)}^2 \leq \lambda_n (x_0;d\mu)^{-1} \int \abs{\pi_n(x)}^2\, d\mu(x)
\end{equation}
or
\begin{equation} \lb{A.1.8}
\abs{\pi_n(x_0)}^p \leq \lambda_n (x_0,p;d\mu)^{-1} \int \abs{\pi_n(x)}^p\, d\mu(x)
\end{equation}

3. Our proof shows the $\inf$ in \eqref{A.1.4} is a $\min$ and the minimizing $\pi$ is
given by
\begin{equation} \lb{A.1.9x}
\pi(x) = K(x,x_0)
\end{equation}
\end{remarks}

\begin{proof} Expand $\pi_n$ in terms of the orthonormal basis $\{p_j\}_{j=0}^n$:
\begin{equation} \lb{A.1.9}
\pi_n(x) =\sum_{j=0}^n a_j p_j(x)
\end{equation}
$\pi_n (x_0)=1$ is equivalent to
\begin{equation} \lb{A.1.10}
\sum_{j=0}^n a_j p_j (x_0) =1
\end{equation}

By the Schwartz inequality
\begin{align}
1 &\leq K(x_0,x_0) \sum_{j=0}^n a_j^2  \lb{A.1.11} \\
&= K(x_0, x_0) \int \pi_n(x)^2\, d\mu(x) \lb{A.1.12}
\end{align}
where equality occurs in \eqref{A.1.11} if $a_j = \ol{p_j(x_0)}/K(x_0,x_0)$, that is,
if $\pi_n$ is given by \eqref{A.1.9}. \eqref{A.1.5} is immediate from this case of
equality and \eqref{A.1.12}.
\end{proof}

Christoffel numbers have been a critical tool in OP theory for over a century, with
important uses by Erd\"os-Tur\'an \cite{ET40} and turned to high art by Freud and
Nevai; see Nevai \cite{NevFr}. They can be used for lower bounds on $K$\!, that is,
upper bounds for $\lambda_n$, by using any convenient trial polynomial for $\pi_n$
(see (d) below). One gets upper bounds on $K$\!, that is, lower bounds for $\lambda_n$,
by the immediate

\begin{corollary}\lb{C.A.2} If $d\mu \geq d\nu$, then
\begin{equation} \lb{A.1.13}
K_n(x_0,x_0;d\mu) \leq K_n (x_0,x_0;d\nu)
\end{equation}
\end{corollary}

\begin{remarks} 1. This shows the true power of Theorem~\ref{T.A.1} and the need
to allow $\mu(\bbR)\neq 1$.

2. In particular, if
\begin{equation} \lb{A.1.14}
d\mu =f(x)\, dx + d\mu_\s
\end{equation}
then
\begin{equation} \lb{A.1.15}
K_n(x_0,x_0;d\mu) \leq K_n (x_0, x_0; f\,dx)
\end{equation}
\end{remarks}

\medskip
\noindent {\bf (b) Gauss-Jacobi quadrature}. \ The main result here is

\begin{theorem}[Gauss-Jacobi quadrature] \lb{T.A.3} Let $\mu$ be an arbitrary
positive nontrivial measure on $\bbR$ with finite moments. Fix $n$ and define $d\mu_n$
to be the point measure with weights only at the zeros $\{x_j^{(n)}\}_{j=1}^n$ of
$p_n(x)$ and weights
\begin{equation} \lb{A.1.16}
d\mu_n (\{x_j^{(n)}\}) = \lambda_n (x_j^{(n)};d\mu)
\end{equation}
the Christoffel numbers of $d\mu$. Then, if $\pi$ is a polynomial of degree $2n-1$
or less, we have
\begin{equation} \lb{A.1.17}
\int \pi(x)\, d\mu(x) =\int \pi(x)\, d\mu_n(x)
\end{equation}
\end{theorem}

\begin{remark} In our applications, we will care much more that the masses are at
the zeros than the variational formulae for the weights.
\end{remark}

\begin{proof}[Sketch of Proof] Here is a proof intended for Schr\"odinger operator
experts. (For the more usual OP proof, see Freud's book \cite{FrB}.) Let $J_{n;F}$ be
an $n\times n$ matrix in the upper corner of the Jacobi matrix, \eqref{1.4}, associated
to $d\mu$. Then the recursion \eqref{1.3} implies that if
\begin{equation} \lb{A.1.18}
u_j (z)=p_{j-1}(z) \qquad j=1,2,\dots, n
\end{equation}
then
\begin{equation} \lb{A.1.19}
[(J_{n;F}-z)u]_j = -a_{n+1} \delta_{jn} p_n(z)
\end{equation}
Thus the eigenvalues of $J_{n;F}$ are the zeros of $p_n$ and the normalized eigenvectors
are $p_{j-1}(z)/K(z,z)^{1/2}$. It follows that $d\mu_n$ is the spectral measure of
$J_{n;F}$ with eigenvector $(1\, 0 \dots 0)^t$.

For any measure $d\eta$ and associated Jacobi matrix, $J$, since $d\eta$ is the spectral
measure for $\delta\equiv (1\, 0 \dots)^t$, we have
\begin{equation} \lb{A.1.20}
\int x^\ell \, d\eta = \langle\delta, J^\ell\delta\rangle
\end{equation}
so $\int x^{\ell+k}\,d\eta = \langle J^k\delta, J^\ell \delta\rangle$ for $k,\ell =
0,\dots, n-1$ depends only on $\{J^k\delta\}_{k=0}^{n-1}$ and so on $\{a_k,b_k\}_{k=0}^{n-1}$.
Moreover,
\begin{equation} \lb{A.1.21}
\int x^{2n-1} \, d\eta =\langle J^{n-1}\delta, J J^{n-1}\delta\rangle
\end{equation}
only depends on $\{a_k,b_k\}_{k=0}^{n-1}\cup\{b_n\}$. Thus, $J$ and $J_{n;F}$, which have the
same set of these parameters, have the same moments of order up to $2n-1$, that is,
\begin{equation} \lb{A.1.22}
\int x^k\, d\mu = \int x^k\, d\mu_{n;F} \qquad 0 \leq k \leq 2n-1
\end{equation}
which is \eqref{A.1.17}.
\end{proof}

\medskip
\noindent {\bf (c) Bernstein's inequality}. \ These inequalities control $\pi'_n$ in
terms of $n$ and $\pi_n$ for polynomials $\pi_n$ of degree at most $n$.

\begin{theorem}\lb{T.A.4} Let $\partial\bbD$ be the unit circle in $\bbC$, $\partial\bbD
=\{z\mid\abs{z}=1\}$. Let $\pi_n$ be a polynomial of degree $n$. Then
\begin{equation} \lb{A.1.23}
\sup_{z\in\partial\bbD}\, \abs{\pi'_n(z)}\leq n\, \sup_{z\in\partial\bbD}\,
\abs{\pi_n(z)}
\end{equation}
\end{theorem}

\begin{remark} If $\pi_n(z) = z^n$, one has equality in \eqref{A.1.23}.
\end{remark}

\begin{proof}(Szeg\H{o} \cite{Sz28}) Since $\pi_n (e^{i\theta})=\sum_{j=0}^n a_j e^{ij\theta}$,
we have
\begin{equation} \lb{A.1.24}
\pi_n (e^{i\theta}) =\int_0^{2\pi} \sum_{j=0}^n e^{ij(\theta-\varphi)} \pi_n
(e^{i\varphi}) \, \f{d\varphi}{2\pi}
\end{equation}
so
\begin{align}
-i\pi'_n (e^{i\theta}) &= \int_0^{2\pi} \sum_{j=1}^n je^{ij(\theta-\varphi)} \pi_n
(e^{i\theta})\, \f{d\rho}{2\pi} \notag \\
&=\int F_n (\theta-\varphi) e^{in(\theta-\varphi)} \pi_n (e^{i\varphi})\,
\f{d\rho}{2\pi} \lb{A.1.25}
\end{align}
where
\begin{equation} \lb{A.1.26}
F_n(\theta) =\sum_{j=-n+1}^{n-1} (n-\abs{j}) e^{ij\theta}
\end{equation}
(for the $j>0$ terms in \eqref{A.1.25}, integrate to zero).

By cancellation,
\begin{equation} \lb{A.1.27}
(1-\cos\theta) F_n(\theta) = 1-\cos (n\theta)
\end{equation}
so $F_n(\theta)\geq 0$ and, by \eqref{A.1.26}, $\int F_n (\theta)\f{d\theta}{2\pi}
=n$. Thus,
\[
\abs{\pi'_n (e^{i\theta})} \leq \|\pi_n\|_\infty \int \abs{F_n(\theta-\varphi)}\,
\f{d\varphi}{2\pi} = n\|\pi_n\|_\infty
\qedhere
\]
\end{proof}

\begin{theorem}\lb{T.A.5} Let $\pi_n$ be an arbitrary polynomial of degree $n$. Then
\begin{equation} \lb{A.1.28}
\sup_{x\in [-a,a]}\, \bigl[\abs{\pi'_n(x)} (a^2-x^2)^{1/2}\bigr] \leq 3n
\sup_{x\in [-a,a]} \, \abs{\pi_n(x)}
\end{equation}
\end{theorem}

\begin{proof} By scaling, we need only check the case $a=2$. Define
\begin{equation} \lb{A.1.29}
\ti\pi_n(z)=z^n \pi_n \biggl( z+\f{1}{z}\biggr)
\end{equation}
$\ti\pi_n$ is a polynomial of degree $2n$ so, by \eqref{A.1.23},
\begin{equation} \lb{A.1.30}
\sup_\theta \, \abs{\ti\pi'_n (e^{i\theta})}\leq 2n \sup_{x\in [-2,2]}\,
\abs{\pi_n(x)}
\end{equation}
since $e^{i\theta}\to e^{i\theta} + e^{-i\theta}=2\cos\theta$ maps $\partial\bbD$
to $[-2,2]$.

By \eqref{A.1.29},
\begin{align*}
\pi'_n\biggl(z+\f{1}{z}\biggr) (1-z^{-2}) &= \f{d}{dz}\, z^{-n}  \ti\pi_n(z) \\
&= -nz^{-n-1} \ti\pi_n(z) + z^{-n} \ti\pi'_n(z)
\end{align*}
so, by \eqref{A.1.30},
\[
\sup_{e^{i\theta}\in\partial\bbD}\, \abs{\pi'_n(2\cos\theta) 2\sin\theta} \leq 3n
\|\ti\pi_n\|_\infty
\]
which is \eqref{A.1.28} for $a=2$.
\end{proof}

\medskip
\noindent {\bf (d) Dirichlet trial polynomials}. \ For use in both \eqref{A.1.4} and
\eqref{A.1.17}, we want a rich set of trial polynomials, $\pi_n(x)$. In particular,
we want $\pi_n$'s concentrated near $x=x_0$ and otherwise small in some interval
$[x_0-a, x_0+a]$. By scaling, we may as well consider $x_0=0$, $a=1$. An analyst
might try $(1-x^2)^n$, but that has width $n^{-1/2}$ --- and we will see that one
can do better. We will get width $n^{-1}$. One can't do better than this, by Bernstein's
inequality, if $\pi_n(\theta)=1$ and $\|\pi_n\|_\infty =1$, then $\pi_n(x) \geq\f12$
for $\abs{x}\leq \f{1}{2n} -O(\f{1}{n^2})$.

Our choice is related to Dirichlet and Fej\'er kernels and is, in fact, essentially
the minimizer for the Christoffel problem with $x_0=0$ and $d\mu=\chi_{[-1,1]}
(1-x^2)^{-1/2}\, dx$.

\begin{theorem}\lb{T.A.6} For any $x_0\in\bbR$ and $a>0$, there exist, for each $n$,
polynomials $\pi_n (x;x_0,a)$ so that
\begin{SL}
\item[{\rm{(i)}}]
\begin{equation} \lb{A.1.31}
 \deg \pi_n =2n-2
\end{equation}

\item[{\rm{(ii)}}]
\begin{equation} \lb{A.1.32}
\pi_n (x_0) =1
\end{equation}

\item[{\rm{(iii)}}]
\begin{equation} \lb{A.1.33}
\abs{\pi_n(x)}\leq \min\biggl(1, \f{1}{2n} + \f{a}{2n\abs{x-x_0}}\biggr)
\qquad \text{if } \abs{x-x_0} \leq a
\end{equation}

\item[{\rm{(iv)}}] For any $\delta \leq a$,
\begin{equation} \lb{A.1.34}
\int_{x_0-\delta}^{x_0+\delta} \abs{\pi_n(x)}^2\, dx = \f{\pi a}{n}
+O\biggl(\f{1}{n^2}\biggr)
\end{equation}
\end{SL}
\end{theorem}

\begin{remark} \eqref{A.1.33} implies $\abs{\pi_n(x)}\leq C_\delta/n$ if $\abs{x-x_0}
>\delta$, and for any $\veps$, $\abs{\pi_n(x)}<\veps$ if $\abs{x-x_0}\geq C_\veps/n$.
\end{remark}

\begin{proof} By scaling, we can suppose that $x_0=0$, $a=1$, in which case
we will call the polynomials $D_n$, that is,
\begin{equation} \lb{A.1.35}
\pi_n (x;x_0,a) =D_n \biggl( \f{x-x_0}{a}\biggr)
\end{equation}
Recall there are polynomials $T_n(x)$ (Chebyshev of the first kind) with
\begin{equation} \lb{A.1.40}
\deg T_n =n
\end{equation}
so that
\begin{equation} \lb{A.1.41}
T_n (\cos\theta) =\cos (n\theta)
\end{equation}
Define $D_n$ by
\begin{equation} \lb{A.1.42}
D_n(x) =\f{1}{n} \, \sum_{j=0}^{n-1} (-1)^j T_{2j}(x)
\end{equation}

By \eqref{A.1.40}, \eqref{A.1.31} holds for $D_n$. By \eqref{A.1.41} (and $\cos\theta
=0\Leftrightarrow \theta =\f{\pi}{2}$ mod $\pi$), $T_{2j}(0)=(-1)^j$, so $D_n$ obeys
\eqref{A.1.32} for $x_0=0$.

By \eqref{A.1.41}, $\abs{T_n(x)}\leq 1$ on $[-1,1]$, so
\[
\abs{D_n(x)}\leq 1 \qquad \text{on } [-1,1]
\]
which is half of \eqref{A.1.33}. For the other half, sum the geometric series to see
that
\begin{equation} \lb{A.1.43}
D_n (\cos\theta) = \f{1}{2n} + \f{(-1)^{n-1}}{2n}\, \f{\cos((2n-1)\theta)}{\cos\theta}
\end{equation}
which implies the other half of \eqref{A.1.33}.

Since
\[
\int_{-\pi}^\pi \cos (2k\theta) \cos (2j\theta) \, \f{d\theta}{2\pi} =
\begin{cases}
1 & \text{if } k=j=0 \\
\f12 & \text{if } k = j\neq 0 \\
0 & \text{if } k\neq j
\end{cases}
\]
we have
\begin{align*}
\int_{-\pi}^\pi D_n^2 (\cos\theta)\, \f{d\theta}{2\pi}
&= \f{1}{n^2}\, \biggl[ 1 + \f12 \, (n-1)\biggr] \\
&= \f{1}{2n} + O\biggl( \f{1}{n^2}\biggr)
\end{align*}
Since $d\theta = (1+O(x^2))\, dx$ near $x=0$ and $D_n^2 = O(\f{1}{n^2})$ away from
$x=0$, we obtain \eqref{A.1.34} when $a=1$.
\end{proof}

\bigskip

\end{document}